\documentclass[final,5p,twocolumn]{elsarticle}

\usepackage{mathrsfs,amsmath,amssymb,amsthm,graphicx,color,dsfont}
\usepackage{multirow,float,algorithm,algpseudocode}
\usepackage{mathtools,cuted,booktabs}

\theoremstyle{definition}

\newtheorem{assumption}{Assumption}
\newtheorem{lemma}{Lemma}
\newtheorem{theorem}{Theorem}
\newtheorem{remark}{Remark}
\newtheorem{problem}{Problem}
\newtheorem{definition}{Definition}

\newtheorem{proposition}{Proposition}

\newcommand{\qedblck}{\hfill\blacksquare}
\renewcommand{\t}{\mbox{\tiny\sf T}}
\newcommand{\R}{\mathbb{R}}

\usepackage{etoolbox}
\makeatletter
\patchcmd{\@makecaption}
  {\scshape}
  {}
  {}
  {}
\makeatother

\begin{document}

\begin{frontmatter}

\title{Data-driven stabilizer design and closed-loop analysis of general nonlinear systems via Taylor's expansion}

\author[add1]{Meichen Guo}
\ead{m.guo@tudelft.nl}
\author[add2]{Claudio De Persis}
\ead{c.de.persis@rug.nl}
\author[add3]{Pietro Tesi}
\ead{pietro.tesi@unifi.it}

\address[add1]{Delft Center for Systems and Control,
Delft University of Technology, 2628 CD Delft, The Netherlands}
\address[add2]{Engineering and Technology Institute Groningen,
University of Groningen, Nijenborgh 4, 9747 AG  Groningen, The Netherlands}
\address[add3]{ Department of Information Engineering, University of Florence, 50139 Florence, Italy}

\begin{abstract}
For data-driven control of nonlinear systems, the basis functions characterizing the dynamics are usually essential. In existing works, the basis functions are often carefully chosen based on pre-knowledge of the dynamics so that the system can be expressed or well-approximated by the basis functions and the experimental data. For a more general setting where explicit information on the basis functions is not available, this paper presents a data-driven approach for stabilizer design and closed-loop analysis via the Lyapunov method. First, based on Taylor's expansion and using input-state data, a stabilizer and a Lyapunov function are designed to render the known equilibrium locally asymptotically stable. Then, data-driven conditions are derived to check whether a given sublevel set of the found Lyapunov function is an invariant subset of the region of attraction. One of the main challenges is how to handle Taylor's remainder in the design of the local stabilizers and the analysis of the closed-loop performance.
\end{abstract}


\end{frontmatter}

\section{Introduction}

Most control approaches of nonlinear systems are based on well-established models of the system constructed by pre-knowledge or system identification. When the models are not explicitly constructed, nonlinear systems can be directly controlled using input-output data. Controlling a system via input-output data without explicitly identifying the model is called the direct data-driven control method, and it has been gaining more and more attentions for both linear and nonlinear systems. An early survey of data-driven control methods can be found in \cite{Hou2013Survey}. More recently, the authors of \cite{tanaskovic2017data} developed an online control approach, the work \cite{Hou2019tac} utilized the dynamic linearization data models for discrete-time non-affine nonlinear systems, the authors of \cite{Tabuada2017CDC} and \cite{Fraile2021feedbacklinearization} considered feedback linearizable systems, and the works \cite{Berberich2022MPCnonlinear} and \cite{Liu2022} designed data-driven model predictive controllers. Inspired by Willems \emph{et al.}'s fundamental lemma, \cite{cdpTAC2020} proposed data-driven control approaches for linear and nonlinear discrete-time systems. Using a matrix Finsler's Lemma, \cite{Waarde2021} applied data-driven control to Lur'e systems. The authors of \cite{Dai2021NonlinearDDcdc} used state-dependent representation and proposed an online optimization method for data-driven stabilization of nonlinear dynamics. For polynomial systems, \cite{PolyTAC} designed global stabilizers using noisy data, and \cite{Nejati2022} synthesized data-driven safety controllers. The recent work \cite{Martin2022polyapprox} investigated dissipativity of nonlinear systems based on polynomial approximation.

\emph{Related works.} Some recent works related to nonlinear data-driven control and the region of attraction (RoA) estimation are discussed in what follows.

Deriving a data-based representation of the dynamics is one of the important steps in data-driven control of unknown nonlinear systems. If the controlled systems are of certain classes, such as polynomial systems having a known degree, the monomials of the state can be chosen as basis functions to design data-driven controllers such as presented in \cite{guoCDC2020,PolyTAC}. By integrating noisy data and side information, \cite{Ahmadi2020} showed that unknown polynomial dynamics can be learned via semidefinite programming. When the nonlinearities satisfy quadratic constraints, data-driven stabilizer was developed in \cite{Luppi2022}. With certain knowledge and assumptions on the nonlinear basis function, systems containing more general types of nonlinearities have also been studied in recent works. For instance, under suitable conditions, some nonlinear systems can be lifted into polynomial systems in an extended state for control, such as the results shown in \cite[Section IV]{Strasser2021beyondpoly} and \cite[Section 3.2]{Kaiser2021}. Using knowledge of the basis functions, \cite{DePersis2022cancellation} designed data-driven controllers by (approximate) cancellation of the nonlinearity. When the system nonlinearities cannot be expressed as a combination of known functions, \cite{DePersis2022cancellation} presented data-driven local stabilization results by choosing basis functions carefully such that the neglected nonlinearities are small in a known set of the state. On the other hand, if the knowledge on the basis functions is not available, approximations of the nonlinear systems are often involved. The previous work \cite{cdpTAC2020} tackled the nonlinear data-driven control problem by linearizing the dynamics around the known equilibrium and obtaining a local stability result. According to these existing results, it is clear that the efficiency and the performance of data-driven controllers can be improved via pre-known knowledge such as specific classes of the systems or the nonlinear basis functions. Nonetheless, there is still a lack of comprehensive investigation of the more general case where the nonlinear basis functions cannot be easily and explicitly obtained. 

The RoA estimation is another relevant topic in nonlinear control. For general nonlinear systems, it is common that the designed controllers only guarantee local stability. Hence, it is of importance to obtain the RoA of the closed-loop systems for the purpose of theoretical analysis as well as engineering applications. Unfortunately, for general nonlinear systems, it is extremely difficult to derive the exact RoA even when the model is explicitly known. A common solution is to estimate the RoA based on Lyapunov functions. Using Taylor's expansion and considering the worst-case remainders, \cite{Chesi09nonPoly} estimated the RoA of uncertain non-polynomial systems via linear matrix inequality (LMI) optimizations. RoA analysis for polynomial systems was presented in \cite{Tan2008tac} using polynomial Lyapunov functions and sum of squares (SOS) optimizations. The authors of \cite{Topcu2009tac} studied uncertain nonlinear systems subject to perturbations in certain forms and used the SOS technique to compute invariant subsets of the RoA. It is noted that, in these works, the RoA estimation winds up in solving bilinear optimization problems, and techniques such as bisection or special bilinear inequality solver are required to find the solutions. For nonlinear systems without explicit models, there are also efforts devoted to learning the RoA by various approaches. The authors of \cite{Bobiti2018tac} developed a sampling-based approach for a class of piecewise continuous nonlinear systems that verifies stability and estimates the RoA using Lyapunov functions. Based on the converse Lyapunov theorem, \cite{Colbert2018cdc} processed system trajectories to lift a Lyapunov function whose level sets lead to an estimation of the RoA. Using the properties of recurrent sets, \cite{Shen2022learningROA} proposed an approach that learns an inner approximation of the RoA via finite-length trajectories. It should be pointed out that, all the aforementioned works focus on stability analysis of autonomous nonlinear systems, \emph{i.e.}, the control design is not considered.

\emph{Contributions.} For general nonlinear systems, this paper presents a data-driven approach to simultaneously obtaining a Lyapunov function and designing a state feedback stabilizer that renders the known equilibrium locally asymptotically stable. Using Taylor's expansion, the unknown dynamics are approximated by linear systems or polynomial systems. Then, linear stabilizers and polynomial stabilizers are designed for the approximated models using finite-length input-state data collected in an off-line experiment. To handle the remainder resulting from the approximation, we conduct the experiment close to the known equilibrium such that the remainder is small with a known bound. An over-approximation of all the feasible dynamics is then found using the collected data, and Petersen's lemma \cite{Petersen1987} is used for the controller design. For polynomial approximations, as the conditions characterizing the stabilizers are positive conditions for polynomial matrices, the SOS technique \cite{Papachristodoulou2005ACC} is utilized to make the conditions easily solvable. The data-driven stabilizer design can be seen as a generalization of the nonlinear control result in the previous work \cite[Section V.B]{cdpTAC2020}. Specifically, this paper considers both the linear and the polynomial approximations of continuous-time systems. Analogous results can also be derived for discrete-time systems. On the other hand, the focus of this paper is on general nonlinear systems, while the previous works \cite{PolyTAC} and \cite{Bisoffi2021petersen} dealt with only polynomial systems. 
In comparison with \cite{DePersis2022cancellation} where the influence of the remaining nonlinearity is attenuated by careful choices of the basis functions, this work achieves the objective by using Taylor's polynomials as basis functions and conduct the experiment close to the equilibrium.

For estimating the RoA with the designed data-driven controller, we first derive an estimation of the remainder by assuming a known bound on the high-order derivatives of the unknown functions. With the help of the Positivstellensatz \cite{Stengle1974Psatz} and the SOS relaxations, we derive data-driven conditions that verify whether a given sublevel set of the obtained Lyapunov function is an invariant subset of the RoA. This is achieved by deriving a sufficient condition for the negativity of the derivative of the Lyapunov function based on the estimation of the remainder. The conditions are derived via data and some prior knowledge on the dynamics, and can be easily solved by software such as Matlab. The estimated RoA gives insights to the closed-loop system under the designed data-driven controller, and is relevant for both theoretical and practical purposes. Note that alternatively, the RoA can be estimated based on other data-driven methods, such as the ones developed in  \cite{Bobiti2018tac,Colbert2018cdc,Shen2022learningROA}. Simulations results on the inverted pendulum show the applicability of the design and estimation approach. For future works, it is of our interests to thoroughly investigate the influence of certain parameters and further improve the design and analysis.


The rest of the paper is arranged as follows. Section \ref{section:preliminaries} formulates the problem and then presents relevant techniques and existing results for the subsequent sections. Data-driven control designs with different orders of approximations are presented in Section \ref{section:controldesign}. The data-driven characterization of the RoA is derived in Section \ref{section:RoAestimation}. Numerical results and analysis on the inverted pendulum are illustrated in Section \ref{section:example}. Finally, Section \ref{section:conclusion} concludes the paper.

\emph{Notation.} Throughout the paper, $A\succ(\succeq)0$ denotes that matrix $A$ is positive (semi-)definite, and $A\prec(\preceq)0$ denotes that matrix $A$ is negative (semi-)definite. For vectors $a,b\in\R^n$, $a\preceq b$ means that $a_i\le b_i$ for all $i=1,\dots,n$. $\|\cdot\|$ denotes the Euclidean norm.

Moreover, we list some important symbols used in subsequent sections in the table below.
\medskip

\begin{tabular}{l l}
  \hline
  $R(x,u)$ & Taylor's remainder of $f(x,u)$\\
  $S=\begin{bmatrix} B & A\end{bmatrix}$ & True dynamics of the first-/high-order\\
  & approximation \\
  $\widehat S=\begin{bmatrix} \widehat B & \widehat A\end{bmatrix}$ & Dynamics consistent with the collect-\\
  & ed data \\
  $\mathcal{C}_k$ & Feasible set of the dynamics at each\\
  & time $t_k$ \\
  $\textbf{A}_k$, $\textbf{B}_k$, $\textbf{C}_k$ & Matrices describing set $\mathcal{C}_k$ \\
  $\mathcal I$ & Intersection of the sets $\mathcal C_k$,\\
  & $k=0,\dots,T-1$\\
  $\overline{\mathcal I}$ & Over-approximation of the set $\mathcal I$ \\
  $\overline{\textbf{A}}$, $\overline{\textbf{B}}$, $\overline{\textbf{C}}$ & Matrices describing set $\overline{\mathcal I}$ \\
  \hline
\end{tabular}

\section{Problem Formulation and Preliminaries}\label{section:preliminaries}

Consider a general nonlinear continuous-time system
\begin{align}\label{dynamics:nonlinear generalCT}
\dot x = f(x,u)
\end{align}
where the state $x\in\R^n$ and the input $u\in\R^m$. The function $f$ is of class $ C^{r}$ for some integer $r\ge 0$. Assume that $(x_e,u_e)$ is a known equilibrium of the system to be stabilized. For simplicity and without loss of generality, in this paper we let $(x_e,u_e)=(0,0)$, as any known equilibrium can be converted to the origin by a change of coordinates.

For a general nonlinear discrete-time system, \cite[Section V.B]{cdpTAC2020} approximates the unknown dynamics as a linear system around a known equilibrium and developed a data-driven local stabilizer. In this work, we show that the result in therein can be extended to approximations with any order via Taylor's expansion. Moreover, for the closed-loop system under the designed data-driven controller, a data-driven estimation of the RoA will also be presented.

To gather information regarding the system, we perform $E$ experiments on the system, where $E$ is an
integer satisfying $1 \le E \le T$ and $T$ an integer equalling the total number of collected samples. On one
extreme, one could perform $1$ single experiment during which a total of $T$ samples are collected. On the
other extreme, one could perform $T$ independent experiments, during each one of which a single sample is
collected. The advantage of short multiple experiments is that they allow the designer to collect information
about the system at different points in the state space without incurring in problems due to the free evolution of the system.

Either way, a dataset $\mathcal{DS}_c:= \{\dot x(t_k); x(t_k); u(t_k)\}^{T-1}_{k=0}$ for the continuous-time system can be obtained. Organize the data collected in the experiment(s) as
\begin{align*}
X_{0} &= \begin{bmatrix}
x(t_{0}) & \cdots & x(t_{T-1})
\end{bmatrix},\\
U_{0} &= \begin{bmatrix}
u(t_0) & \cdots & u(t_{T-1})
\end{bmatrix},\\
X_{1} &= \begin{bmatrix}
\dot x(t_{0}) & \cdots & \dot x(t_{T-1})
\end{bmatrix}.
\end{align*}

The problem studied in this work is formulated as follows.

\begin{problem}[\emph{Data-driven stabilizer design and RoA estimation}]
For system \eqref{dynamics:nonlinear generalCT}, design a state feedback controller $u=k(x)$ using the input-state data $X_0$, $X_1$, and $U_0$, such that the origin is locally asymptotically stable for the closed-loop system, and an inner estimation of the RoA of the closed-loop system is derived.
\end{problem}

\begin{remark}[\emph{On the experimental data}]
The derivative data $\dot x(t_k)$, $k=1,\dots,T-1$ can be approximated using numerical differentiation. For instance, using the forward difference approximation gives
\[
\dot x_i(t_k) = \frac{x_i(t_{k+1})-x_i(t_{k})}{t_{k+1}-t_k} +e_i(t_k), ~~i=1,\dots,N,
\]
where the approximation error $e_i(t_k)$ is proportional to $t_{k+1}-t_k$ and can be handled as noise.

In this work, the data is assumed to be noiseless for the sake of simplicity. If the data is corrupted by noise, to decrease the effect of noise in the derivative approximation, numerical approaches such as the total variation regularization \cite{Rudin1992totalvariation} can be used as shown in \cite{Chartrand2011} and \cite{Brunton2016}.
 $\qedblck$
\end{remark}

The rest of this section presents techniques and existing results needed for the subsequent data-driven control design and analysis.

\subsection{Taylor's expansion}\label{subsec:Taylor}
To approximate the unknown nonlinear dynamics, we use Taylor's expansion to represent the unknown function by the sum of the $r$th Taylor polynomial and the remainder for any integer $r\ge 0$ at a given equilibrium. As explained previously, we assume that the functions in consideration have an equilibrium at the origin.

Consider a function $\phi:\R^\sigma\rightarrow \R$ that is of class $C^{r}$ on an open convex set $\mathbb D\subseteq \R^n$ containing the origin, and $\phi(0)=0$. Define $\alpha=(\alpha_{1},\dots,\alpha_{\sigma})$ as a $\sigma$-tuple of nonnegative integers and
\begin{align*}
|\alpha| &= \alpha_{1}+\cdots+\alpha_{\sigma},\\
\alpha ! &= \alpha_1 ! \cdots \alpha_\sigma !,\\
~~\partial^{\alpha}\phi &= \frac{\partial^{|\alpha|}\phi}{\partial z_1^{\alpha_{1}}\cdots\partial z_\sigma^{\alpha_{\sigma}}},\\
z^\alpha &= z_1^{\alpha_1} \cdots z_\sigma^{\alpha_\sigma}
\end{align*}
for any $z=[z_1~\cdots~z_{\sigma}]^{\top}\in\R^{\sigma}$. As given in \cite{FollandTaylor2}, for $z\in\mathbb D$, Taylor's expansion of the function $\phi$ at the origin is
\begin{align}
\phi(z) &= \sum_{|\alpha|\le r}\frac{\partial^{\alpha}\phi(0)}{\alpha !}z^{\alpha} + R_r(z)
\end{align}
for any integer $r\ge 0$, and the remainder can be expressed in integral form as
\begin{align}\label{TaylorRemainder}
R_r(z)=(r+1)\sum_{|\alpha|=r+1}\frac{z^{\alpha}}{\alpha !}\int^{1}_{0} (1-t)^r \partial^{\alpha}\phi(tz) dt.
\end{align}
In \cite{FollandTaylorRemainder}, the author has shown that the estimate of the modulus of $\partial \phi^{\alpha}(z)$, $|\alpha|=r$ at the origin can be used to obtain an estimate for the remainder $R_r(z)$.

\begin{lemma}\label{lemma:remainderbound}
Consider all $z\in\mathbb{D}$ where $\mathbb D\subseteq\mathbb R^\sigma$ is an open convex set containing the origin. If for $|\alpha|=r$ with any integer $r\ge 0$, $\phi$ is of class $C^{r}$ for all $z\in\mathbb{D}$, and there exists $C\ge 0$ such that
\begin{align}\label{continuouscondition:scalar}
|\partial^{\alpha}\phi(z) - \partial^{\alpha}\phi(0)| \le C\|z\|~~\forall z\in\mathbb{D},
\end{align}
then the remainder $R_r(z)$ in \eqref{TaylorRemainder} satisfies
\begin{align}
\left| R_r(z) \right| \le \frac{\sigma^{r/2}C\|z\|^{r+1}}{(r+1)!}~~\forall z\in\mathbb{D}.
\end{align}
\end{lemma}

The proof of Lemma \ref{lemma:remainderbound} can be found in \ref{proofrmndrbnd}.

\subsection{Petersen's lemma}
In the section for data-driven controller design, Petersen's lemma is essential for deriving the sufficient condition characterizing the controller. Due to the space limit, the proof of the lemma is omitted\ and one may refer to works such as \cite{Petersen1987,Shcherbakov2008PetersenLemma,Bisoffi2021petersen} for more details.

\begin{lemma}[Petersen's lemma \cite{Petersen1987}]
\label{lemma:Petersen}
Consider matrices $\mathcal G=\mathcal G^{\top}\in\R^{n\times n}$, $\mathcal M\in\R^{n\times m}$, $\mathcal M\ne 0$, $\mathcal N \in\R^{p\times n}$, $\mathcal N \ne 0$, and a set $F$ defined as
\[
F =\{\mathcal F\in\R^{m\times p}: \mathcal F^{\top} \mathcal F\preceq \overline{\mathcal F} \}
\]
where $\overline{\mathcal F}=\overline{\mathcal F}^{\top}\succeq 0$. Then, for all $\mathcal F\in F$
\[
\mathcal G + \mathcal M \mathcal F \mathcal N + \mathcal N^{\top} \mathcal F^{\top} \mathcal M^{\top} \preceq 0
\]
holds if and only if there exists $\mu>0$ such that
\[
\mathcal G + \mu \mathcal M \mathcal M^{\top} +\mu^{-1}\mathcal N^{\top}\overline{\mathcal F}\mathcal N \preceq 0.
\]
\end{lemma}

\subsection{Sum of squares relaxation}
As solving positive conditions of multivariable polynomials is in general NP-hard, the SOS relaxations are often used to obtain sufficient conditions that are tractable. The SOS polynomial matrices are defined as follows.

\begin{definition}
(\emph{SOS polynomial matrix} \cite{Chesi2010Survey})
$M:\R^n\rightarrow\R^{\sigma\times \sigma}$ is an SOS polynomial matrix if there exist $M_{1},\dots,M_{k}:\R^n\rightarrow\R^{\sigma\times \sigma}$ such that
\begin{align}\label{def_SOSmatrix}
  M(x)=\sum^{k}_{i=1}M_{i}(x)^{\top}M_{i}(x)~~\forall x\in\R^n.
\end{align}
\end{definition}

Note that when $\sigma =1$, $M(x)$ becomes a scalar SOS polynomial.

It is straightforward to see that if a matrix $M(x)$ is an SOS polynomial matrix, then it is positive semi-definite, \emph{i.e.},
$M(x)\succeq 0~\forall x\in\R^n$.
Relaxing the positive polynomial conditions into SOS polynomial conditions makes the conditions tractable and easily solvable by common software.

\subsection{Positivstellensatz}

In the RoA analysis, we need to characterize polynomials that are positive on a semialgebraic set and the Positivstellensatz plays an important role in the characterization.

Let $p_1,\dots,p_k$ be polynomials. The multiplicative monoid, denoted by $\mathcal S_M(p_1,\dots,p_k)$, is the set generated by taking finite products of the polynomials $p_1,\dots,p_k$. The cone $\mathcal S_C(p_1,\dots,p_k)$ generated by the polynomials is defined as
\begin{align*}
&\mathcal S_C(p_1,\dots,p_k)\\
& =\{ s_0 +\sum_{i=1}^j s_iq_i:s_0,\dots,s_j \text{ are SOS polynomials},\\
&\quad\quad\quad\quad\quad\quad\quad\quad q_1,.\dots,q_j\in\mathcal S_M(p_1,\dots,p_k) \}.
\end{align*}
The ideal $\mathcal S_I(p_1,\dots,p_k)$ generated by the polynomials is defined as
\[
\mathcal S_I(p_1,\dots,p_k)= \left\{ \sum_{i=1}^k r_ip_i: r_1,\dots,r_k \text{ are polynomials} \right\}.
\]
Stengle's Positivstellensatz \cite{Stengle1974Psatz} is presented as follows in \cite{Chesi2010Survey}.

\begin{theorem}[\emph{Positivstellensatz}]\label{theorem:P-satz}
Let $f_1,\dots,f_k$, $g_1,\dots,g_l$, and $h_1,\dots,h_m$ be polynomials. Define the set
\begin{align*}
\mathcal X = \{x\in\R^n: &f_1(x)\ge 0, \dots, f_k(x)\ge 0, \\
&g_1(x)= 0, \dots, g_l(x)= 0,\\
&\text{and } h_1(x)\ne 0, \dots, h_m(x)\ne 0 \}.
\end{align*}
Then, $\mathcal X = \emptyset$ if and only if
\[
\exists f\!\in\mathcal S_C(f_1,\dots,f_k),g\!\in\mathcal S_I(g_1,\dots,g_l),h\!\in\mathcal S_M(h_1,\dots,h_m)
\]
such that
\[
f(x)+g(x)+h(x)^2=0.
\]
\end{theorem}

For the subsequent RoA analysis, we will use the following result derived from the Positivstellensatz.

\begin{lemma}\label{lemma:Psatz}
Let $\varphi_1$ and $\varphi_2$ be polynomials in $x$. If there exist SOS polynomials $s_1$ and $s_2$ in $x$ such that \begin{align}\label{soscondition:Psatzlemma}
-(s_1\varphi_1(x)+s_2\varphi_2(x) +x^{\t}x) \text{ is SOS }~\forall x\in\R^n
\end{align}
then the set inclusion condition
\begin{align}\label{setcondtion:RoA}
\{x\in\R^n: \varphi_1(x)\ge0, x\ne 0\} \subseteq \{ x\in\R^n: \varphi_2(x)<0 \}
\end{align}
holds.
\end{lemma}

The proof of Lemma \ref{lemma:remainderbound} can be found in \ref{proofPsatz}.

\section{Data-driven controller design}\label{section:controldesign}

To approximate the system \eqref{dynamics:nonlinear generalCT}, we write Taylor's expansion of the function $f(x,u)$ as
\begin{align*}
f(x,u)=f_{appr}(x,u)+R(x,u),
\end{align*}
where $f_{appr}(x,u)$ contains Taylor's polynomials up to a certain degree and $R(x,u)$ represents the truncated remainder containing high-order terms. This remainder constitutes the main uncertainty of the data-based representation of the closed-loop system and brings difficulties in the controller design. To attenuate the influence of the remainder, we collect state data close to the equilibrium such that the maximum amplitude of the remainder during the experiment is instantaneously bounded by a known constant. While the first-order approximation of the system leads to a solution to Problem 1, to further minimizing the remainder in the neighborhood of the equilibrium, we can truncate Taylor's polynomials at a higher degree and obtain a high-order approximation.

Using the collected data, we first find an over-approximation of the set containing all dynamics that are consistent with the data. Then, stabilizers are designed such that the origin is locally asymptotically stable for all the dynamics in the approximated set.

\subsection{First-order approximation}

First, for completeness, we will address the local data-driven stabilizer design via first-order approximation.

Consider the continuous-time nonlinear system \eqref{dynamics:nonlinear generalCT}. Denote each element of $f$ as $f_i$ and let $f_i \in C^1(\R^n\times \R^m)$, $i=1,\dots,n$.
The first-order approximation of \eqref{dynamics:nonlinear generalCT} through Taylor's expansion of $f(x,u)$ is
\begin{align}\label{system:approxCT}
\dot x = Ax +Bu +R(x,u)
\end{align}
where $R(x,u)$ denotes the remainder and
\begin{align*}
A = \left. \frac{\partial f(x,u)}{\partial x}\right|_{(x,u)=(0,0)},~~B=\left. \frac{\partial f(x,u)}{\partial u}\right|_{(x,u)=(0,0)}.
\end{align*}
One can treat the remainder as a disturbance that affects the data-driven characterization of the unknown dynamics, and focus on controlling the approximated linear dynamics to obtain a local linear stabilizer. Then, by tackling the impact of $R(x,u)$, the RoA of the closed-loop system can be also be characterized.

%

Similar to our previous works such as \cite{cdpTAC2020} and \cite{PolyTAC}, some bound on the experimental data of the remainder is needed to design a controller. For this purpose, the following assumption is posed.

\begin{assumption}\label{assumption:instanRemainder}
For $k=0,\dots, T$ and a known $\gamma$,
\begin{align}\label{remainderbound:Assumption1}
R(x(t_k),u(t_k))^{\top}R(x(t_k),u(t_k)) \le \gamma^2.
\end{align}
\end{assumption}

\begin{remark}[\emph{Instantaneous remainder bound}]
Assumption \ref{assumption:instanRemainder} is an instantaneous bound on the maximum amplitude of the remainder during the experiment. The bound can be obtained by prior knowledge of the model, such as the physics of the system. If such knowledge is unavailable, one may resort to an over-estimation of $\gamma$. 
$\qedblck$
\end{remark}

\begin{remark}[\emph{Noisy data.}] If the data is corrupted by additive measurement and/or actuator noise, \emph{i.e.}, $x(t_k) = x^*(t_k) +d_x(t_k)$ and $u(t_k) = u^*(t_k) +d_u(t_k)$ where $x^*(t_k)$, $u^*(t_k)$ are the true data, and $d_x(t_k)$, $d_u(t_k)$ represent the noise, the proposed design and analysis approach is still applicable to the noisy case, provided that the remainder $R(x(t_k),u(t_k))$ in Assumption 1 is replaced by $R(x^*(t_k),u(t_k))+d(t_k)$, where $d(t_k)$ is the total noise due to $d_x(t_k)$, $d_u(t_k)$, and the derivative approximation error. 
 $\qedblck$
\end{remark}

\subsubsection{Over-approximation of the feasible set}

Under Assumption \ref{assumption:instanRemainder}, an over-approximation of the set of dynamics that are consistent with the experimental data can be derived as shown in \cite{Bisoffi2021tradeoffs}.

Denote $S=[B~~A]$ as the true dynamics. Based on \eqref{system:approxCT}, at each time $t_{k}$, $k=0,\dots,T-1$, the collected data satisfies
\begin{align*}
\dot x(t_k) = S \begin{bmatrix} u(t_k) \\  x(t_k)\end{bmatrix} +R\big( x(t_k),u(t_k) \big).
\end{align*}
Under Assumption \ref{assumption:instanRemainder}, one has that
\begin{align*}
&\quad R\big( x(t_k),u(t_k) \big)R\big( x(t_k),u(t_k) \big)^{\top}\\
 &= \left( \dot x(t_k) - S \begin{bmatrix} u(t_k) \\  x(t_k)\end{bmatrix} \right)\left( \dot x(t_k) - S \begin{bmatrix} u(t_k) \\  x(t_k)\end{bmatrix} \right)^{\top} \preceq \gamma^2 I.
\end{align*}
Hence, at each time $t_{k}$, $k=0,\dots,T-1$, the matrices $\widehat S=[\widehat B~~\widehat A]$ consistent with the data belongs to the set
\begin{align}
\mathcal C_k = \left\{ \widehat S: \textbf{C}_k+\widehat S\textbf{B}_k+\textbf{B}_k^{\top}\widehat S^{\top}+\widehat S\textbf{A}_k\widehat S^{\top} \preceq 0 \right\}
\end{align}
where
\begin{align*}
\textbf{A}_k &= l(t_k) l(t_k)^{\top},~~ \textbf{B}_k=-l(t_k)\dot x(t_k)^{\top}, \\
\textbf{C}_k &=\dot x(t_k)\dot x(t_k)^{\top}-\gamma^2 I,~~  l(t_k)=\begin{bmatrix} u(t_k) \\ x(t_k)\end{bmatrix}.
\end{align*}
Then, the feasible set of matrices $\widehat S$ that is consistent with all data collected in the experiment(s) is the intersection of all the sets $\mathcal C_k$, \emph{i.e.}, $\mathcal I = \bigcap^{T-1}_{k=0}\mathcal C_k$. 
Though the exact set $\mathcal I$ is difficult to obtain, an over-approximation of $\mathcal I$ in the form of a matrix ellipsoid and of minimum size can be computed. Denote the over-approximation set as
\begin{align}\label{setbarI}
\overline {\mathcal I}:=\left\{ \widehat S: \overline{\textbf{C}}+\widehat S\overline{\textbf{B}}+\overline{\textbf{B}}^{\top}\widehat S^{\top}+\widehat S\overline{\textbf{A}}\widehat S^{\top} \preceq 0 \right\}
\end{align}
where $\overline{\textbf A}=\overline{\textbf A}^{\top}\succ 0$, $\overline{\textbf C}$ is set as $\overline{\textbf C}=\overline{\textbf B}^{\top}\overline{\textbf A}^{-1}\overline{\textbf B}^{\top}-\delta I$ and $\delta>0$ is a constant fixed arbitrarily. Following \cite[Section 5.1]{Bisoffi2021tradeoffs}, the set $\overline{\mathcal I}$ can be found by solving the optimization problem
\begin{align}\label{fsblsetOvrapprox}
\underset{\overline{\textbf A}, \overline{\textbf B}, \overline{\textbf C}}{\text{minimize}} \quad & -\mathrm{log~det}(\overline{\textbf A})\\
\text{subject to}\quad  & \overline{\textbf A}=\overline{\textbf A}^{\top}\succ 0\notag \\
& \tau_k\ge 0,~k=0,\dots,T-1 \notag\\
&\begin{bmatrix}
- \delta I -\sum\limits^{T-1}_{k=0} \tau_k \textbf C_k & \overline{\textbf B}^{\top}-\sum\limits^{T-1}_{k=0} \tau_k \textbf B_k^{\top} & \overline{\textbf B}^{\top}\\
\overline{\textbf B}-\sum\limits^{T-1}_{k=0} \tau_k \textbf B_k & \overline{\textbf A}-\sum\limits^{T-1}_{k=0} \tau_k \textbf A_k & 0\\
\overline{\textbf B} & 0 & -\overline{\textbf A}
\end{bmatrix} \preceq 0\notag.
\end{align}

\begin{remark}[\emph{On the parameter $\delta$}]
To make sure that the over-approximation set $\overline{\mathcal{I}}$ is not empty, we need that $\overline{\textbf B}^{\top}\overline{\textbf A}^{-1}\overline{\textbf B}^{\top}-\overline{\textbf C}\succ0$, which is guaranteed by setting $\overline{\textbf C}=\overline{\textbf B}^{\top}\overline{\textbf A}^{-1}\overline{\textbf B}^{\top}-\delta I$ for some positive constant $\delta$. As pointed out in \cite[Section 3.7]{Boyd1994LMI}, the description of the set $\overline{\mathcal{I}}$ is homogeneous, that is, the matrices $\overline{\textbf A}$, $\overline{\textbf B}$, and $\overline{\textbf C}$ can be scaled by any positive factor without affecting $\overline{\mathcal{I}}$. Hence, the positive constant $\delta$ can be fixed arbitrarily without changing the resulting $\overline{\mathcal{I}}$. In many works, such as \cite{Boyd1994LMI} and \cite{Bisoffi2021tradeoffs}, the variables are normalized, \emph{i.e.}, $\delta$ is set as $1$. In this work, we keep the parameter $\delta$ because although it does not affect the over-approximation of the feasible set, it has effects on the subsequent RoA estimation as shown in the simulation results. However, it is still unclear to us how $\delta$ affects the RoA estimation quantitatively or how to choose an optimal $\delta$ for the design and analysis. These questions can only be answered after more thorough studies on the parameters in the control design process and the RoA analysis, which is out of the scope of this work. $\qedblck$
\end{remark}

\begin{remark}[\emph{Persistency of excitation}]
\label{rmk:PE}
As pointed out in \cite[Section 3.1]{Bisoffi2021tradeoffs}, if the collected data is rich enough, \emph{i.e.}, $\begin{bmatrix} U_0 \\ X_0 \end{bmatrix}$ has full row rank, then the intersection set $\mathcal I$ is bounded, which allows the optimization problem \eqref{fsblsetOvrapprox} to have a solution. Hence, $\begin{bmatrix} U_0 \\ X_0 \end{bmatrix}$ having full row rank implies the feasibility of \eqref{fsblsetOvrapprox}.
\end{remark}


\subsubsection{Data-driven stabilizer design }

Stabilizing the linear approximation of the unknown system renders the origin locally asymptotically stable as the remainder $R(x,u)$ contains higher-order terms and converges to the origin faster than the linear part in a neighborhood of the origin. Hence, the objective of the controller design is to stabilize the origin for all dynamics belonging to the over-approximation set $\overline{\mathcal I}$. This can be achieved in the same manner as done in \cite{Bisoffi2021petersen} via Petersen's lemma. For the completeness of this work, we give the following result on designing a data-driven local stabilizer. The proof shares the same idea with Theorem 2 in \cite{Bisoffi2021petersen}, and thus is omitted in this paper.

\begin{theorem}\label{thrm:controldesignCT}
Under Assumption \ref{assumption:instanRemainder}, given a constant $w>0$, if there exist matrices $Y$ and $P=P^{\top}$ such that
\begin{equation}
\begin{aligned}
\begin{bmatrix}
wP-\overline{\textbf{C}} & \overline{\textbf{B}}^{\top}- \begin{bmatrix} Y \\ P \end{bmatrix}^{\top}\\
\overline{\textbf{B}}- \begin{bmatrix} Y \\ P \end{bmatrix} &  -\overline{\textbf{A}}
\end{bmatrix}&\preceq 0\\
P&\succ 0,\label{LMI:thrmCntrlCT}
\end{aligned}
\end{equation}
then the origin is a locally asymptotically stable equilibrium for the closed-loop system composed of \eqref{dynamics:nonlinear generalCT} and the control law $u=YP^{-1}x$.
\end{theorem}

Consider the Lyapunov function $V(x)=x^{\top}P^{-1}x$. The controller designed by Theorem \ref{thrm:controldesignCT} guarantees that the derivative of $V(x)$ along the trajectory of the closed-loop system $\dot x =\widehat A x+ \widehat B u$ satisfies that
\[
\dot V(x)\le -wV(x)
\]
for any given constant $w>0$ and any $\begin{bmatrix} \widehat B & \widehat A \end{bmatrix}\in\overline{\mathcal I}$. Hence, by choosing the value of $w$, a certain decay rate of the closed-loop solution is enforced.

\subsection{High-order approximation}

One may also approximate the nonlinear system as a polynomial system by truncating Taylor's expansion at a degree higher than $1$. This will result in a smaller remainder in the neighborhood close to the origin. Similar to the first-order approximation case, by making sure that the remainder converges to the origin faster, a stabilizer for the approximated system can render the origin locally asymptotically stable for the overall system.

To write the approximated system into a linear-like form, we consider the nonlinear input-affine system
\begin{align}\label{dynamics:nonlinear affine CT}
\dot x = f(x)+g(x)u
\end{align}
where $x\in\R^n$, $u\in\R^m$, and $f(0)=0$. Functions $f$ and $g$ are of class $C^{r_f}$ and $C^{r_g}$ respectively for some integers $r_f,r_g\ge1$. Using Taylor's expansion, we truncate the series of polynomials at order $r_f$ for functions $f_i$, and at order $r_g$ for functions $g_{ij}$, $i=1,\dots,n$, $j=1,\dots,m$, respectively. Following Section \ref{subsec:Taylor}, the functions can be written as
\begin{align*}
f_i(x) &= \sum_{|\alpha|\le r_f}\frac{\partial^{\alpha}f_i(0)}{\alpha !}x^{\alpha} + R_{f_i}(x),\\
g_{ij}(x) &= \sum_{|\alpha|\le r_g}\frac{\partial^{\alpha}g_{ij}(0)}{\alpha !}x^{\alpha} + R_{g_{ij}}(x)
\end{align*}
with remainders $R_{f_i}(x)$ and $R_{g_{ij}}(x)$.
Then, one can write the polynomial part into the linear-like form and obtain the system
\begin{align}\label{system:linearlike}
\dot x &=  AZ(x) +BW(x)u +R(x,u)
\end{align}
where $R(x,u)= [R_{1}(x,u) \dots R_{n}(x,u)]^{\top}$,
\begin{align}\label{remainder:highorder}
R_i(x,u) = R_{f_i}(x) + \sum^{m}_{j=1}R_{g_{ij}}(x)u_j,
\end{align}
$Z(x)$ is a vector of monomials in $x$ having degree $1$ to $r_f$, $W(x)$ is a matrix of monomials in $z$ having degree $0$ to $r_g$, constant matrices $A$, $B$ are unknown. Denote the degree of $u$ as $r_u$, then the remainder $R_i(x,u)$ has the degree $r_R:=\max\{r_f,r_g+r_u\}$.

\subsubsection{Over-approximation of the feasible set}

Similar to the case of first-order approximation, when Assumption \ref{assumption:instanRemainder} holds, one can obtain an over-approximation of the feasible set by solving the optimization problem \eqref{fsblsetOvrapprox}. Abusing the notations $\overline{\textbf{A}}$, $\overline{\textbf{B}}$, and $\overline{\textbf{C}}$, we define the over-approximation set as
\begin{align}\label{tildebarI}
\widetilde {\mathcal I}:=\left\{ \widehat S: \overline{\textbf{C}}+\widehat S\overline{\textbf{B}}+\overline{\textbf{B}}^{\top}\widehat S^{\top}+\widehat S\overline{\textbf{A}}\widehat S^{\top} \preceq 0 \right\},
\end{align}
where $\overline{\textbf C}=\overline{\textbf B}^{\top}\overline{\textbf A}^{-1}\overline{\textbf B}^{\top}-\delta I$ for some arbitrarily fixed constant $\delta>0$, $\overline{\textbf{A}}$ and $\overline{\textbf{B}}$ are the solutions to the optimization problem \eqref{fsblsetOvrapprox}, with
\begin{align*}
\textbf{A}_k &= \ell(t_k) \ell(t_k)^{\top},~~ \textbf{B}_k=-\ell(t_k)\dot x(t_k)^{\top}, \\
\textbf{C}_k &=\dot x(t_k)\dot x(t_k)^{\top}-\gamma^2 I,~~  \ell(t_k)=\begin{bmatrix} W\big(x(t_k) \big)u(t_k) \\ Z\big(x(t_k)\big)\end{bmatrix}.
\end{align*}

Similar to the analysis in Remark \ref{rmk:PE}, if the data is rich enough, \emph{i.e.}, $\begin{bmatrix} \overline U_0 \\ Z_0 \end{bmatrix}$ has full row rank, then the set $\mathcal{I}$ is bounded, allowing \eqref{fsblsetOvrapprox} to have a solution.

Note that in the linear like forms \eqref{system:linearlike}, the vector $Z(x)$ contains all monomials in $x$ having degree from $1$ to $r_f$. Hence, the size of $Z(x)$ can be substantially large due to high truncating degree or high system order. As pointed out in our previous work \cite{PolyTAC}, a large vector $Z(x)$ tends to cause computational issues in solving the sufficient conditions characterizing the stabilizers. Therefore, to avoid bringing a large $Z(x)$ directly into the conditions, we use a smaller vector $\widehat Z(x)$ for the Lyapunov function and controller design, same as in \cite{PolyTAC}. Specifically, $\widehat Z(x)$ is chosen as a $p\times 1$ vector on $x$ such that $\widehat Z(x)=0$ if and only if $x=0$ and $Z(x)=H(x)\widehat Z(x)$ where the matrix $H(x)$ is non-unique. Then, the linear-like system is rewritten as
\begin{align*}
\dot x &=  AH(x)\widehat Z(x) +BW(x)u +R(x,u).
\end{align*}
Define the Lyapunov function as $V(x)= \widehat Z(x)^{\top}P^{-1}\widehat Z(x)$ where $P\succ 0$ and the controller as $u=K(x)\widehat Z(x)$. Note that the vector $\widehat Z(x)$ needs to be chosen such that $V(x)$ is radially unbounded, and a simple choice of $\widehat Z(x)$ is to make its first $n$ components coincide with $x$.

\begin{table*}[t]
\begin{align}
&\Upsilon(x)= -\frac{\partial\widehat Z(x)}{\partial x} \left(-\overline{\textbf A}^{-1}\overline{\textbf B}\right)^{\top}
\begin{bmatrix}  W(x)Y(x) \\ H(x)P  \end{bmatrix}
 -\! \begin{bmatrix}  W(x)Y(x) \\ H(x)P  \end{bmatrix}^{\top} \! \left(-\overline{\textbf A}^{-1}\overline{\textbf B}\right)\frac{\partial\widehat Z(x)}{\partial x}^{\top} \!
-\mu(x)\delta\frac{\partial\widehat Z(x)}{\partial x}\frac{\partial\widehat Z(x)}{\partial x}^{\top}
\label{Upsilon_CT}\\
&\frac{\partial\widehat Z(x)}{\partial x} \left( -\overline{\textbf{A}}^{-1}\overline{\textbf{B}} \right)^{\top}\begin{bmatrix}  W(x)Y(x) \\ H(x)P  \end{bmatrix}
+(\star)^{\top}+\epsilon(x)I
+\mu(x) \delta \frac{\partial\widehat Z(x)}{\partial x} \frac{\partial\widehat Z(x)}{\partial x}^{\top}
\!+\mu(x)^{-1}\!\begin{bmatrix}  W(x)Y(x) \\ H(x)P  \end{bmatrix}^{\top} \!\overline{\textbf{A}}^{-1}\!
\begin{bmatrix}  W(x)Y(x) \\ H(x)P  \end{bmatrix} \preceq 0 \label{MItheoremHighorderProof}
\end{align}
\end{table*}

\subsubsection{Data-driven stabilizer design}

Based-on the over-approximated set $\widetilde{\mathcal{I}}$, we can characterize sufficient conditions for a local stabilizer using Lyapunov's second method and Petersen's lemma. To make the condition tractable, the SOS technique is applied to relax the positivity conditions for the polynomial matrices. The design of the data-driven stabilizer is given in the following result.

\begin{theorem}\label{thrm:cntrldsgnCTHighorder}
Under Assumption \ref{assumption:instanRemainder}, if there exist polynomial $\mu(x)>0$ $\forall x\in\R^n$, a positive definite polynomial $\epsilon(x)$, \emph{i.e.}, zero at the origin and positive elsewhere, and matrices $Y(x)$ and $P$, such that
\begin{align}\label{SOScondition:cntrlCT}
\begin{bmatrix}
\Upsilon(x) -\epsilon(x)I_p & \begin{bmatrix}  W(x)Y(x) \\ H(x)P  \end{bmatrix}^{\top}\\
\begin{bmatrix}  W(x)Y(x) \\ H(x)P  \end{bmatrix} & \mu(x)\overline{\textbf A}
\end{bmatrix} \text{  is SOS}
\end{align}
where $\Upsilon(x)$ is defined in \eqref{Upsilon_CT}
with the parameter $\delta$ used to define $\overline{\textbf C}$ for the set $\widetilde{\mathcal I}$, then the state feedback controller $u=Y(x)P^{-1}\widehat Z(x)$ makes the origin a locally asymptotically stable equilibrium for the system $\dot x = f(x)+g(x)u$.
\end{theorem}

\proof
First, following the description \eqref{tildebarI} of set $\widetilde{\mathcal I}$ and the definition of $\overline{\textbf{C}}$, it holds that, for all $\widehat S\in\widetilde{\mathcal I}$,
\begin{align*}
&\quad ~~ \overline{\textbf{C}}+\widehat S\overline{\textbf{B}}+\overline{\textbf{B}}^{\top}\widehat S^{\top}+\widehat S\overline{\textbf{A}}\widehat S^{\top}\\
&=\left( \widehat S^{\top}+\overline{\textbf{A}}^{-1}\overline{\textbf{B}} \right)^{\top}\overline{\textbf{A}} \left( \widehat S^{\top}+ \overline{\textbf{A}}^{-1}\overline{\textbf{B}} \right) -\overline{\textbf{B}}^{\top}\overline{\textbf{A}}^{-1}\overline{\textbf{B}}+\overline{\textbf{C}}\\
&=\left( \widehat S^{\top}+\overline{\textbf{A}}^{-1}\overline{\textbf{B}} \right)^{\top}\overline{\textbf{A}} \left( \widehat S^{\top}+ \overline{\textbf{A}}^{-1}\overline{\textbf{B}} \right) -\delta I \preceq 0.
\end{align*}
Define $\Delta = \overline{\textbf{A}}^{1/2}\left( \widehat S^{\top}+ \overline{\textbf{A}}^{-1}\overline{\textbf{B}} \right)$, and it follows that $\Delta^{\top}\Delta \preceq \delta I$ and
\[
\widehat S^{\top} = -\overline{\textbf{A}}^{-1}\overline{\textbf{B}} + \overline{\textbf{A}}^{-1/2}\Delta.
\]

The objective is to find a control gain $K(x)=Y(x)P^{-1}$ that stabilizes the approximated systems for all $\widehat S\in \widetilde{\mathcal I}$. Hence, the closed-loop system of the controlled approximated system is
\begin{align*}
\dot x&= \widehat AZ(x)+\widehat BW(x)Y(x)P^{-1}\widehat Z(x) \\
&= \widehat S \begin{bmatrix}  W(x)Y(x) \\ H(x)P  \end{bmatrix} P^{-1}\widehat Z(x).
\end{align*}
Taking the time derivative of the Lyapunov function $V(x)= \widehat Z(x)^{\top}P^{-1}\widehat Z(x)$ along the trajectory of the closed-loop approximated system gives
\begin{align*}
\dot V(x) = 
 \widehat Z(x)^{\top}P^{-1}
\frac{\partial\widehat Z(x)}{\partial x} \widehat S
\begin{bmatrix}  W(x)Y(x) \\ H(x)P  \end{bmatrix} P^{-1}\widehat Z(x) +(\star)^{\top}.
\end{align*}
Recall that $\widehat S^{\top} = -\overline{\textbf{A}}^{-1}\overline{\textbf{B}} + \overline{\textbf{A}}^{-1/2}\Delta$. Let $\Phi(x)$ be such that $\dot V(x) =\widehat Z(x)^{\top}P^{-1}\Phi(x) P^{-1} \widehat Z(x)$. Then, it holds true that
\begin{align*}
\Phi(x) &= \frac{\partial\widehat Z(x)}{\partial x} \left( -\overline{\textbf{A}}^{-1}\overline{\textbf{B}} \right)^{\top}\begin{bmatrix}  W(x)Y(x) \\ H(x)P  \end{bmatrix}\\
& \quad +\frac{\partial\widehat Z(x)}{\partial x} \left( \overline{\textbf{A}}^{-1/2}\Delta \right)^{\top}
\begin{bmatrix}  W(x)Y(x) \\ H(x)P  \end{bmatrix}  +(\star)^{\top}.
\end{align*}

If the SOS condition in Theorem \ref{thrm:cntrldsgnCTHighorder} is satisfied, using the Schur complement, one derives \eqref{MItheoremHighorderProof}.
By Petersen's Lemma (Lemma \ref{lemma:Petersen}), the inequality \eqref{MItheoremHighorderProof} is equivalent to $\Phi(x)\preceq -\epsilon(x)I_p$.

For the original systems with remainders, one has that
\begin{align*}
\dot V(x) & \le -\epsilon(x)\widehat Z(x)^{\top}P^{-2}\widehat Z(x)\\
&\quad +2\widehat Z(x)^{\top}P^{-1}\frac{\partial\widehat Z(x)}{\partial x} R\big( x,Y(x)P^{-1}\widehat Z(x) \big).
\end{align*}
Note that the remainder $R\big( x,Y(x)P^{-1}\widehat Z(x) \big)$ having the form \eqref{remainder:highorder} is of degree $r_R=\max\{r_f,r_g+r_u\}$. As the degrees of $\epsilon(x)$ and $\widehat Z(x)$ are fixed by design, we can guarantee that the term $2\widehat Z(x)^{\top}P^{-1}\frac{\partial\widehat Z(x)}{\partial x} R\big( x,Y(x)P^{-1}\widehat Z(x) \big)$ is of higher order and thus converges to $0$ faster than $-\epsilon(x)\widehat Z(x)^{\top}P^{-2}\widehat Z(x)$ for all $x$ in a neighborhood of the origin. Moreover, as $\epsilon(x)>0$ $\forall x\ne 0$, the origin is locally asymptotically stable for the closed-loop system $\dot x = f(x)+g(x)Y(x)P^{-1}\widehat Z(x)$. \qed

\begin{remark}[\emph{Dimension of $\widehat Z(x)$}]
The feasibility of the condition $\Phi(x)\preceq -\epsilon(x)I_p$ when $x\ne 0$ implies that rank$(\Phi(x))= p$ for $x\ne 0$. Meanwhile, since $\widehat Z(x)=0$ if and only if $x=0$, one has that $p\ge n$. Hence,  rank$(\frac{\partial \widehat Z(x)}{\partial x})\le n$, and it holds that rank$(\Phi(x))\le 2n$ $\forall x\in\R^n$. Therefore, the feasibility of $\Phi(x)\preceq -\epsilon(x)I_p$ implies that the dimension $p$ of vector $\widehat Z(x)$ is such that $n\le p\le 2n$. $\qedblck$
\end{remark}

\begin{table*}[t]
{\small
\begin{align}
\label{inequality:HODT}
\begin{bmatrix}
-P+\epsilon(x)I_n & \begin{bmatrix}  W(x)Y(x) \\ H(x)P  \end{bmatrix}^{\top} \left(-\overline{\textbf A}^{-1}\overline{\textbf B}\right) \\
\star & -P
\end{bmatrix} +\mu(x)\delta  \begin{bmatrix} 0 \\ I_n \end{bmatrix} \begin{bmatrix} 0 \\ I_n \end{bmatrix}^{\top}
+\mu(x)^{-1} \begin{bmatrix}
\begin{bmatrix}  W(x)Y(x) \\ H(x)P  \end{bmatrix}^{\top} \overline{\textbf A}^{-1/2} \\ 0
\end{bmatrix}
\begin{bmatrix}
\begin{bmatrix}  W(x)Y(x) \\ H(x)P  \end{bmatrix}^{\top} \overline{\textbf A}^{-1/2} \\ 0
\end{bmatrix}^{\top}\preceq 0 
\end{align}
}
\end{table*}

\begin{remark}[\emph{Comparison of Theorem \ref{thrm:cntrldsgnCTHighorder} with previous results}]
Theorem \ref{thrm:cntrldsgnCTHighorder} provides a local data-driven controller design for general nonlinear dynamics. It is a continuous-time generalization of the previous result \cite[Theorem 6]{cdpTAC2020} where the first-order approximation and a linear controller is considered. Compared to results that relies on specific choices of the basis functions, such as \cite{DePersis2022cancellation}, Theorem \ref{thrm:cntrldsgnCTHighorder} synthesizes a data-driven controller using Taylor polynomials as basis functions. The approach for the control of polynomial systems in Theorem \ref{thrm:cntrldsgnCTHighorder} is an alternative to the results presented in \cite{PolyTAC} and \cite{Bisoffi2021petersen}, where \cite{PolyTAC} handles the additional noisy term in a different way, and \cite{Bisoffi2021petersen} searches for a Lyapunov function without restricting to a special form.  $\qedblck$
\end{remark}

\section{Region of attraction estimation}\label{section:RoAestimation}

In the previous section, we have shown that data-driven stabilizers can be designed for nonlinear systems using first-order or high-order approximations. The resulting controllers make the origin locally asymptotically stable. Besides this property, it is of paramount importance to estimate the RoA of the closed-loop system. The definition of the RoA is given as follows.

\begin{definition}
[\emph{Region of attraction}]
For the systems $\dot x =f(x)$, if for every initial condition $x(t_0)\in\mathcal R$, it holds that $\lim_{t\rightarrow\infty}x(t)=0$, then $\mathcal R$ is a region of attraction of the system with respect to the origin. If there exists a $C^{1}$ function $V:\R^{n}\rightarrow\R$ and a positive constant $c$ such that
\begin{align*}
\Omega_{c}:=\{ x\in\R^{n}:V(x)\le c \}
\end{align*}
is bounded and
\begin{align*}
&V(0)=0,~~V(x)>0~~\forall x\in\R^{n}\\
&\{x\in\R^{n}:V(x)\le c, x\ne0\} \subseteq \{ x\in\R^{n}: \dot V(x)<0 \},
\end{align*}
then $\Omega_{c}$ is an invariant subset, or called an estimation, of the RoA.
\end{definition}

In this section, for the designed data-driven controllers in Section \ref{section:controldesign}, we derive data-driven conditions to determine whether a given sublevel set of the Lyapunov function is an invariant subset of the RoA. The derived conditions are data-driven because they are obtained using the over-approximated set $\overline{\mathcal{I}}$. We note that once the controller is computed, it is possible to use any other data-driven method to estimate the RoA, see for example \cite{Bobiti2018tac,Colbert2018cdc,Shen2022learningROA}.

\subsection{First-order approximation}


By the controller design method in Section \ref{section:controldesign}, the Lyapunov function $V(x)$, and thus the set $\Omega_c$, are available for the analysis. To characterize the set $\{ x\in\R^n : \dot V(x) < 0 \}$, we need a bound on the remainder for $x$ in a neighborhood of the origin. This is achievable by posing the following assumption on the partial derivative of each $f_i$.

\begin{assumption}\label{assumption:Lipschitzr1}
For all $z\in \mathbb D\subseteq \R^{n+m}$ where $\mathbb D$ is a star-convex neighborhood of the origin, $f_i$ is continuously differentiable and
\begin{align}\label{condition:Lipschitzr1}
\left| \frac{\partial f_i}{\partial z_j}(z)-\frac{\partial f_i}{\partial z_j}(0) \right| \le L_i \|z\|~~ \forall j=1,\dots,m+n,
\end{align}
for $i=1,\dots,n$ with known $L_i>0$.
\end{assumption}

\begin{remark}[\emph{Existence of $L_i$}]
Assumption \ref{assumption:Lipschitzr1} is the weakest condition needed for deriving a bound on the remainder using Lemma \ref{lemma:remainderbound}. A stronger condition, such as the Lipschitz continuity of $\frac{\partial f_i}{\partial z_j}$, guarantees the existence of $L_i$. It is also noted that $L_i$ can be estimated using a data-based bisection procedure as shown in \cite[Section III.C]{Martin2022polyapprox}.
$\qedblck$
\end{remark}


Under Assumption \ref{assumption:Lipschitzr1}, it follows from Lemma \ref{lemma:remainderbound} that the first order approximation remainder $R(x,u)$ of $f(x,u)$ satisfies that, for all $(x,u)\in\mathbb D$
\begin{align}
|R_i(x,u)| \le \frac{\sqrt{m+n} L_i}{2} \|(x,u)\|^{2},~~i=1,\dots,n,
\end{align}
where $R_i(x,u)$ is the $i$th element of vector $R(x,u)$.


\begin{remark}[\emph{On Assumptions \ref{assumption:instanRemainder} and \ref{assumption:Lipschitzr1}}]
Under Assumption \ref{assumption:Lipschitzr1}, using Lemma \ref{lemma:remainderbound}, the bound $\gamma^2$ in Assumption \ref{assumption:instanRemainder} can be derived for the experimental data. During the experiment(s), suppose that the smallest ball containing $(x(t),u(t))$ has radius $R_e$, \emph{i.e.}, $\|(x(t_k),u(t_k))\|\le R_e$ for all $k=0,\dots,T$. Then, for $k=0,\dots,T$,
\begin{align*}
&\quad~ R(x(t_k),u(t_k))^{\top}R(x(t_k),u(t_k)) \\
&= \sum^n_{i=1}R_i(x(t_k),u(t_k))^2\\
&\le \sum^n_{i=1} \frac{(m+n)L_i^2}{4}\|x(t_k),u(t_k)\|^4\\
&\le \sum^n_{i=1} \frac{(m+n)L_i^2}{4}R_e^4.
\end{align*}
Hence, if some prior knowledge on the dynamics is known such that Assumption \ref{assumption:Lipschitzr1} holds, $\gamma^2$ can be chosen as $\sum^n_{i=1} \frac{(m+n)L_i^2}{4}R_e^4$ to satisfy Assumption \ref{assumption:instanRemainder}. $\qedblck$

\end{remark}

After finding an estimate of the remainder $R(x,u)$, we analyze the RoA of the closed-loop system.

Under Assumption \ref{assumption:Lipschitzr1}, for the closed-loop system with the controller $u=Kx$ designed via Theorem \ref{thrm:controldesignCT}, the derivative of the Lyapunov function can be described in the following lemma.

\begin{lemma}\label{lemma:dVCT1st}
Consider system \eqref{dynamics:nonlinear generalCT} and the linear controller $u=YP^{-1}x$ where $Y$ and $P$ are designed to satisfy \eqref{LMI:thrmCntrlCT} with any given constant $w>0$. Under Assumption \ref{assumption:Lipschitzr1}, the derivative of the Lyapunov function $V(x)=x^{\top}P^{-1}x$ satisfies that, for all $x\in\mathbb D$,
\begin{align}\label{bound:dVx1st}
\dot V(x) \le -wx^{\top}P^{-1}x + 2\kappa(x)\rho(x)
\end{align}
where
\begin{align}\label{definition:kappa}
\kappa(x):=\begin{bmatrix} x^{\top}Q_1\|(x,Kx)\|^{2} & \cdots & x^{\top}Q_n\|(x,Kx)\|^{2} \end{bmatrix},
\end{align}
$Q_i(x)$ is the $i$th column of $P^{-1}$, and the vector $\rho(x)$ is contained in the polytope
\begin{align}\label{definition:polytopeHCT}
\mathcal H:= \{ \varrho:  -\bar h \preceq \varrho \preceq \bar h \}
\end{align}
with
\[
\bar h=\begin{bmatrix}\bar h_1~ \cdots ~\bar h_n\end{bmatrix}^{\top}
 = \begin{bmatrix}
 \frac{\sqrt{m+n} L_1}{2} &\cdots & \frac{\sqrt{m+n} L_n}{2}
\end{bmatrix}^{\top}.
\]
\end{lemma}

The proof of Lemma \ref{lemma:dVCT1st} can be found in \ref{proofdVCT1st}.

Denote the number of distinct vertices of $\mathcal H$ as $\nu$ and each vertex of $\mathcal H$ as $ h_k$, $k=1,\dots,\nu$. Using the Positivstellensatz result in Lemma \ref{lemma:Psatz}, we present the following result.

\begin{proposition}\label{prpstn:ROAestimationCT}
Suppose that the controller $u=Kx$ renders the origin a locally asymptotically stable equilibrium for \eqref{dynamics:nonlinear generalCT} with the Lyapunov function $V(x)=x^{\top}P^{-1}x$. Under Assumption \ref{assumption:Lipschitzr1}, given a $c>0$ such that $\Omega_c=\{ x\in\R^{n}:V(x)\le c \}\subseteq \mathbb D$, if there exist SOS polynomials $s_{1k},s_{2k}$ in $x$, $k=1,\dots,\nu$ such that
\begin{align}\label{SOSpolyPropCT}
-\left[s_{1k}(c-V(x))+s_{2k}\left(-wx^{\top}P^{-1}x  + 2 \kappa(x) h_k \right)+x^{\top}x \right]
\end{align}
is SOS, where $\kappa(x)$ is defined as in \eqref{definition:kappa} and $h_k$ are the distinct vertices of the polytope $\mathcal H$ defined in \eqref{definition:polytopeHCT}, then $\Omega_c$ is an invariant subset of the RoA of the system $\dot x = f(x,Kx)$ relative to the equilibrium $x = 0$.
\end{proposition}

\proof According to \cite[Page 87]{Blanchini2008SetControl}, the polytope $\mathcal H$ can be expressed as
\begin{align*}
\mathcal H = \left\{ \varrho=\sum^{\nu}_{k=1}\lambda_k(x) h_k,~\sum^{\nu}_{k=1}\lambda_k(x)=1,~ \lambda_k(x)\ge 0 \right\}
\end{align*}
for any fixed $x\in\Omega_c$.
Then, the derivative of the Lyapunov function satisfies that
\begin{align*}
\dot V(x) &\le -wx^{\top}P^{-1}x + 2\kappa(x)\rho(x) \\
&= -wx^{\top}P^{-1}x +2\kappa(x)\sum^{\nu}_{k=1}\lambda_k(x)h_k\\
&= \sum^{\nu}_{k=1}\lambda_k(x)(-wx^{\top}P^{-1}x)+ \sum^{\nu}_{k=1}\lambda_k(x) 2\kappa(x)h_k\\
&= \sum^{\nu}_{k=1}\lambda_k (x)\left(  -wx^{\top}P^{-1}x+2\kappa(x)h_k \right)
\end{align*}
As $\lambda_k(x)\ge 0$ and $\sum^{\nu}_{k=1}\lambda_k(x)=1$, if
\[
-wx^{\top}P^{-1}x+2\kappa(x)h_k<0
\]
holds for all $k=1,\dots,\nu$, then $\dot V(x)<0$.

By Lemma \ref{lemma:Psatz}, for each $k=1,\dots,\nu$, if there exist SOS polynomials $s_{1k},s_{2k}$ such that
\begin{align*}
-\left[s_{1k}(c-V(x))+s_{2k}\left(-wx^{\top}P^{-1}x  + 2 \kappa(x)h_k \right)+x^{\top}x\right]
\end{align*}
is SOS, then the set inclusion condition
\begin{align*}
\{x\in\R^n : &V(x)\le c,~x\ne 0\} \\
&\subseteq \{ x\in\R^n : -wx^{\top}P^{-1}x  + 2 \kappa(x)h_k < 0 \}
\end{align*}
holds. This leads to the set inclusion condition
\begin{align*}
\{x\in\R^n :  V(x)\le c,~x\ne 0\}\subseteq \{ x\in\R^n : \dot V(x) < 0 \},
\end{align*}
and hence $\Omega_c$ is an inner estimate of the ROA.  \qed
\medskip

\begin{remark}[\emph{Numerical method for RoA estimation}]
In \cite{DePersis2022cancellation}, the RoA is estimated by a numerical method, \emph{i.e.}, a sufficient condition of $V(x^+)-V(x)<0$ is found using data, and the grids in a compact region are tested to see whether the sufficient condition is satisfied so that the sublevel sets of $V(x)$ can be found as the RoA estimation. In this work, we derive SOS conditions for the RoA estimation that is an alternative to the mesh method used in \cite{DePersis2022cancellation}. 
\end{remark}


\subsection{High-order approximation}
For high-order approximation, to analyze the closed-loop system, we also need a bound on the high-order remainder $R_f(x)+R_g(x)u$. To simplify the analysis, this subsection considers system \eqref{dynamics:nonlinear affine CT} with a single input, that is $u\in\R$. In this case, the remainder takes the form
\[
R(x,u) = R_f(x) +R_g(x)u
\]
where
\begin{align*}
R(x,u)&=[R_{1}(x,u)~\cdots~R_{n}(x,u)]^{\top},\\
 R_f(x)&=[R_{f_1}(x)~\cdots~R_{f_n}(x)]^{\top}, \text{ and} \\
 R_g(x)&=[R_{g_1}(x)~\cdots~ R_{g_n}(x)]^{\top}.
\end{align*}

To obtain a bound on the remainder using Lemma \ref{lemma:remainderbound}, we pose an assumption on the high-order partial derivatives of the functions $f(x)$ and $g(x)$. 

\begin{assumption}\label{assumption:Continuityr>1affine}
For all $x\in \mathbb D\subseteq \R^{n}$ where $\mathbb D$ is a star-convex neighborhood of the origin, $\partial^{\alpha} f_i$, $\partial^{\beta} g_{i}$, $|\alpha|=r_f$, $|\beta| =r_g$, are absolutely continuous for $x\in\mathbb D$ and
\begin{align*}
|\partial^{\alpha} f_i(x)-\partial^{\alpha} f_i(0)|&\le L_i \|x\|,\\
|\partial^{\beta} g_{i}(x)-\partial^{\beta} g_{i}(0)|&\le M_{i} \|x\|~~i=1,\dots,n,
\end{align*}
with known $L_i,M_{i}>0$.
\end{assumption}

Under Assumption \ref{assumption:Continuityr>1affine}, by Lemma \ref{lemma:remainderbound}, it holds that for all $x\in\mathbb D$
\begin{align}
|R_{f_i}(x)| &\le \frac{\sqrt nL_i}{(r_f+1)!}\|x\|^{r_f+1}, \\
|R_{g_i}(x)| &\le \frac{\sqrt nM_{i}}{(r_g+1)!}\|x\|^{r_g+1}.
\end{align}


Denote the degree of the control input $u$ as $r_u$ and note that there is no constant terms in $u$. Hence, it holds that
\begin{align*}
|u| \le \sum^{r_u}_{j=1} \overline K_{j}\|x\|^{j}
\end{align*}
for some positive constants $\overline K_j$ obtained from the designed $u=K(x)\widehat Z(x)$.
The term $R_{g_i}(x)u$ can be bounded by
\begin{align*}
|R_{g_i}(x)u| &\le  |R_{g_i}(x)||u| \\
&\le \frac{\sqrt n M_i}{(r_g+1)!}\|x\|^{r_g+1} \cdot \sum^{r_u}_{j=1} \overline K_{j}\|x\|^{j} \\
& = \frac{\sqrt n M_i}{(r_g+1)!}\sum^{r_u}_{j=1}\overline K_{j} \|x\|^{r_g+1+j}.
\end{align*}
For any $r_g\ge 1$ and $r_u\ge 1$, one can always find a number $q\ge 1$ and positive constants $\widetilde K_j$ such that
\begin{align}
\sum^{r_u}_{j=1}\overline K_{j} \|x\|^{r_g+1+j} \le \sum^{q}_{j=1} \widetilde K_j \|x\|^{2j}\le \widetilde K_{M}\sum^{q}_{j=1} \|x\|^{2j}
\end{align}
where $\widetilde K_{M} = \max_{j}{\widetilde K_j}$. (Note: for any odd degree term $\|x\|^{a+b}$, we can use the triangular inequality to bound it as $\|x\|^{a+b}\le \epsilon \|x\|^{2a}+\epsilon^{-1}\|x\|^{2b}$ $\forall \epsilon>0$.) We then bound the remainder $R_i(x,K(x)\widehat Z(x))$ as
\begin{align*}
&\quad~ |R_i(x,K(x)\widehat Z(x))|\\
 &\le |R_{f_i}(x)|+|R_{g_i}(x)u| \\
  &\le \frac{\sqrt n L_i}{(r_f+1)!}\|x\|^{r_f+1} +\frac{\sqrt n M_i\widetilde K_{M}}{(r_g+1)!}\sum^{q}_{j=1} \|x\|^{2j}.
\end{align*}
To make sure the obtained bound is a polynomial, we need to set $r_f$ as an odd number. Letting
\[
\widetilde\varrho_i = \max \left\{ \frac{\sqrt n L_i}{(r_f+1)!}, \frac{\sqrt n M_i\widetilde K_{M}}{(r_g+1)!} \right\},
\]
one can write
\[
R_i(x,K(x)\widehat Z(x)) = \widetilde\rho_i(x)\left( \|x\|^{r_f+1} +\sum^{q}_{j=1} \|x\|^{2j} \right)
\]
where $|\widetilde\rho_i(x)|\le \widetilde\varrho_i$ $\forall x\in\mathbb D$.

Denote $Q_i(x)$ as the $i$th column of $P^{-1}\frac{\partial \widehat Z(x)}{\partial x}$, and one has the following term in the derivative of the Lyapunov function
\begin{align*}
&\quad~ 2\widehat Z(x)^{\top}P^{-1}\frac{\partial\widehat Z(x)}{\partial x} R\big( x,K(x)\widehat Z(x) \big)\\
&=2\begin{bmatrix}  \widehat Z(x)^{\top}Q_1(x) & \cdots & \widehat Z(x)^{\top}Q_n(x) \end{bmatrix} R\big( x,K(x)\widehat Z(x) \big)\\
& =2 \sum^n_{i=1} \widehat Z(x)^{\top}Q_i(x)R_i\big( x,K(x)\widehat Z(x) \big).
\end{align*}
Defining
\begin{align}\label{definition:kappaHighorder}
 \widetilde\kappa(x) = \Big[ &\widehat Z(x)^{\top}Q_1(x)( \|x\|^{r_f+1} +\sum^{q}_{j=1} \|x\|^{2j})
~ \cdots \notag\\
&\quad\widehat Z(x)^{\top}Q_n(x)( \|x\|^{r_f+1} +\sum^{q}_{j=1} \|x\|^{2j}) \Big]
\end{align}
and $\widetilde\rho(x) = \begin{bmatrix} \widetilde\rho_1(x) & \cdots & \widetilde\rho_n(x) \end{bmatrix}^{\top}$ gives
\begin{align*}
2\widehat Z(x)^{\top}P^{-1}\frac{\partial\widehat Z(x)}{\partial x} R\big( x,Y(x)P^{-1}\widehat Z(x) \big)
= 2\widetilde\kappa(x)\widetilde\rho(x),
\end{align*}
and thus
\begin{align}\label{bound:dVxhigh}
\dot V(x) \le -\epsilon(x)\widehat Z(x)^{\top}P^{-2}\widehat Z(x)+2\widetilde\kappa(x)\widetilde\rho(x)
\end{align}
for all $x\in\mathbb D$, where the vector $\widetilde\rho(x)$ is contained in the polytope
\begin{align}\label{definition:polytopeHHighorder}
\widetilde{\mathcal H} := \left[-\widetilde\varrho_1,\widetilde\varrho_1 \right]\times\cdots\times
\left[-\widetilde\varrho_n, \widetilde\varrho_n \right].
\end{align}

Following the same spirit of Proposition \ref{prpstn:ROAestimationCT}, we obtain the following result.

\begin{proposition}\label{proposition:ROAHighorder}
Suppose that the $u=Y(x)P^{-1}\widehat Z(x)$ renders the origin a locally asymptotically stable equilibrium for \eqref{dynamics:nonlinear affine CT} with the Lyapunov function $V(x)=\widehat Z(x)^{\top}P^{-1}\widehat Z(x)$. Under Assumption 2, given $c>0$ such that $\Omega_c=\{ x\in\R^{n}:V(x)\le c \}\subseteq\mathbb D$, if there exist SOS polynomials $s_{1k},s_{2k}$ in $x$, $k=1,\dots,\nu$ such that,
\begin{align}\label{SOSpolyPropHighorder}
-\left[s_{1k}(c-V(x))+s_{2k}(-\epsilon(x)\widehat Z(x)^{\top}P^{-2}\widehat Z(x)\right.\notag\\
\left. + 2\widetilde \kappa(x)\widetilde h_k )+x^{\top}x\right] \text{ is SOS}
\end{align}
where $\widetilde\kappa(x)$ is defined as in \eqref{definition:kappaHighorder} and $\widetilde h_k$ are the distinct vertices of polytope $\widetilde{\mathcal H}$ defined in \eqref{definition:polytopeHHighorder}, then $\Omega_c$ is an invariant subset of the RoA of the system $\dot x = f(x)+g(x)Y(x)P^{-1}\widehat Z(x)$ relative to the equilibrium $x = 0$.
\end{proposition}

We summarize the data-driven approach for controller design and RoA estimation in Algorithm \ref{alg:ddcontrol}.

\begin{algorithm}
\caption{Data-driven controller design and RoA estimation}\label{alg:ddcontrol}
\begin{algorithmic}
\item[1.] Collect data $X_0$, $X_1$, and $U_0$ from the offline open-loop experiment(s)
\item[2.] Solve the optimization problem \eqref{fsblsetOvrapprox} for $\overline{\textbf A}$, $\overline{\textbf B}$, and $\overline{\textbf C}$ to obtain $\overline{\mathcal I}$
\item[3.] Solve the LMI \eqref{LMI:thrmCntrlCT} in Theorem 2 ( the SOS condition \eqref{SOScondition:cntrlCT} in Theorem 3) for $Y$ ($Y(x)$) and $P$ to obtain the
   controller $u=YP^{-1}x$ ($u=Y(x)P^{-1}\widehat Z(x)$)
\item[4.] Bound $\dot V(x)$ as in \eqref{bound:dVx1st} (\eqref{bound:dVxhigh}) for the closed-loop system using the designed $u$ and $V(x)$
\item[5.] Determine if $\Omega_c=\{ x\in\R^{n}:V(x)\le c \}$ with a given $c>0$ belongs to the RoA by checking the SOS condition \eqref{SOSpolyPropCT} (\eqref{SOSpolyPropHighorder}) for the first-order (high-order) approximation
\end{algorithmic}
\end{algorithm}


\section{Numerical example}\label{section:example}

In this section, the proposed control design and closed-loop system analysis is applied to the nonlinear benchmark, an inverted pendulum. Simulation results are presented to illustrate the main result of this work.

\subsection{First-order approximation}
Consider the inverted pendulum having the dynamics
\begin{align}\label{dynamics:InvertedPendulum}
\dot x_1 &= x_2,\notag\\
\dot x_2 &= \frac{mgl}{J}\mathrm{sin}(x_1) -\frac{r}{J}x_2 +\frac{l}{J}\mathrm{cos}(x_1) u.
\end{align}
The system is in the form of a general nonlinear system with $f(x,u)=[f_1(x,u)~~f_2(x,u)]^{\top}$ where
\[
f_1(x,u) = x_2,~~f_2(x,u) = \frac{mgl}{J}\mathrm{sin}(x_1) -\frac{r}{J}x_2+\frac{l}{J}\mathrm{cos}(x_1) u.
\]

To collect data, an experiment is conducted with $x(0)= [0.01~~-0.01]^{\top}$ and $u=0.1\mathrm{sin}(t)$ during the time interval $[0,5]$. The data is sampled with fixed sampling period $T_s=0.05$s. The system parameters are $m=0.1$, $g=9.8$, $r=1$, $l=1$, and $J=1$. Collect the data and arrange it into data sets with length $T=10$. We assume that the remainder is over-approximated by $100\%$; in other words, the bound $\gamma$ is twice the largest instantaneous norm of the remainder during the experiment. Then, for the experimental data, Assumption \ref{assumption:instanRemainder} holds with
\[
\gamma = 3.3352\cdot 10^{-6}.
\]
Setting $\delta=0.01$, we first solve the optimization problem \eqref{fsblsetOvrapprox} to find $\overline{\mathcal I}$, and then apply Theorem \ref{thrm:controldesignCT} with $w=1$. The solution found by CVX is
\begin{align}\label{cntrl:exampleCT}
P &= 10^3\cdot \begin{bmatrix}   1.0152  &   -1.3289 \\   -1.3289  &   1.7727 \end{bmatrix},\notag \\
u &= -12.0432x_1  -8.887x_2.
\end{align}

The value of the constant is obtained as follows. We focus on $f_2$ whose first order partial derivatives are
\begin{align*}
\frac{\partial f_2}{\partial x_1} &= \frac{mgl}{J}\mathrm{cos}(x_1)-\frac{l}{J}\mathrm{sin}(x_1)u,\\
\frac{\partial f_2}{\partial x_2} &= -\frac{r}{J},~
\frac{\partial f_2}{\partial u} =\frac{l}{J}\mathrm{cos}(x_1).
\end{align*}
Assumption \ref{assumption:Lipschitzr1} holds for $f_2$ with
\[
L_2 = \sqrt 2 ~\max\left( \frac{mgl}{J},\frac{l}{J} \right).
\]
Assume that $L_2$ is over-estimated by $20\%$, which gives the estimated bound on the remainder $R_2(x,u)$ as
\[
|R_2(x,u)| \le 1.2\cdot \frac{\sqrt{3} L_2}{2}\|(x,u)\|^2=1.4697\|(x,u)\|^2.
\]

Applying Proposition \ref{prpstn:ROAestimationCT}, the largest $c$ found for the controller \eqref{cntrl:exampleCT} is $c^*=7.58\cdot 10^{-4}$.

\begin{figure}
  \centering
  \includegraphics[width=0.4\textwidth]{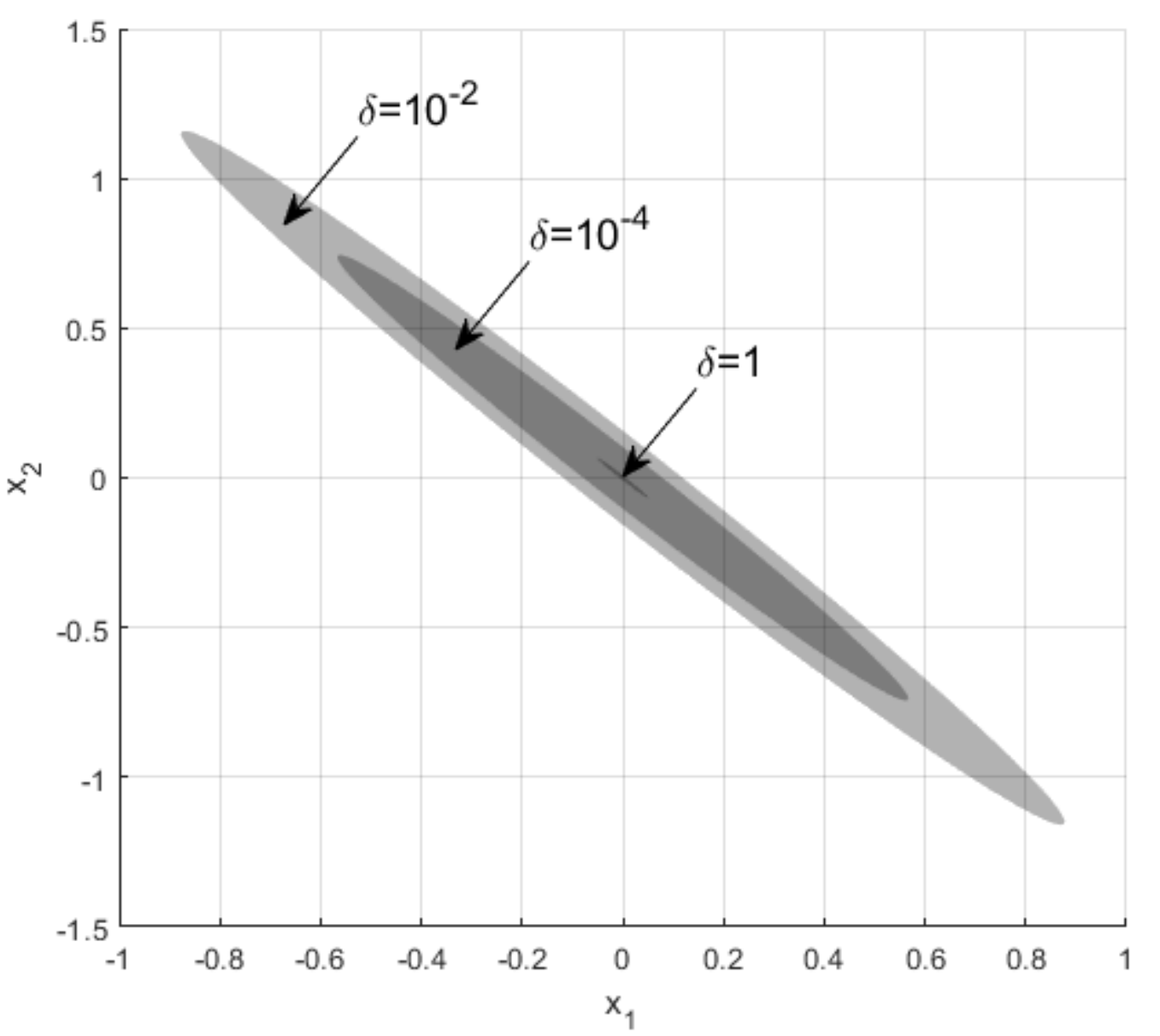}\\
  \caption{The sets $\Omega_c$ found by Proposition \ref{prpstn:ROAestimationCT} with different values of $\delta$. }\label{Fig:ROA1stCTdeltavary}
\end{figure}

To see the effect brought by the parameter $\delta$ in the over-approximation of the feasible set, the same design and analysis approach is repeated for different values of $\delta$. The simulation results associated with various $\delta$ are recorded in Table \ref{table:CTvarydelta}. The resulting sets $\Omega_c$ are illustrated in Figure \ref{Fig:ROA1stCTdeltavary}. As observed from the simulation results, the control gains $K$ do not show significant changes with different values of $\delta$, while the sizes of the estimated RoA vary. The role of $\delta$ in the RoA estimation remains unclear at the moment and interesting questions arise from the simulation results. For instance, how $\delta$ affects the RoA estimation and how to find an optimal $\delta$ for the RoA estimation. It is of interest to answer these questions via a careful and thorough study in the future. 

\begin{table}
\small
\begin{tabular}{l  c  c  c}
  \toprule
   $\delta$ & $P$ & $K$ & $c^*$\\
    \toprule
  $1$ & $10^3\cdot \!\begin{bmatrix} 3.13 & -4.08 \\-4.08  &  5.43\end{bmatrix} $ & $ [-11.84~-8.75] $ & $0.90\cdot \! 10^{-6}$\\
  \hline
  $10^{-2}$ & $ 10^3\cdot\! \begin{bmatrix}1.02 & -1.33\\-1.33 &  1.77\end{bmatrix} $ & $[-12.04~-8.89]$ &$7.58\cdot \! 10^{-4}$\\
  \hline
  $10^{-4}$ & $ 10^3\cdot\! \begin{bmatrix}0.427 &-0.55\\ -0.55 & 0.74\end{bmatrix} $ & $ [-11.67~-8.63] $ & $7.54\cdot \! 10^{-4}$\\
  \bottomrule
\end{tabular}
 \caption{Simulation results on the continuous-time model with different values of $\delta$.}
 \label{table:CTvarydelta}
\end{table}

\subsection{High-order approximation}
Consider again the inverted pendulum \eqref{dynamics:InvertedPendulum}. Setting degrees of the Taylor polynomials as $r_f=5$ and $r_g=2$, one can write the dynamics as
\begin{align*}
\dot x_1 &= x_2,\\
\dot x_2 &=\frac{mgl}{J}x_1 -\frac{mgl}{6J}x_1^3 +\frac{mgl}{120J}x_1^5 -\frac{r}{J}x_2 \\
&\quad +\left( \frac{l}{J}-\frac{l}{2J} x_1^2 \right) u + R_{f_2}(x)+R_{g_2}(x)u
\end{align*}
where the remainders are
\begin{align*}
R_{f_2}(x) &= \left. \frac{mgl}{J}\cdot \sin(\sigma x_1)\right|_{\sigma\in(0,1)}\cdot  \frac{x_1^6}{6!},\\
R_{g_2}(x) &=  \left. \frac{l}{J}\cdot  \cos(\sigma x_1)\right|_{\sigma\in(0,1)}\cdot  \frac{x_1^3}{3!}.
\end{align*}
We again assume that the remainder data is over-approximated by $100\%$, and hence Assumption \ref{assumption:instanRemainder} holds with
\[
\gamma = 2.1602\cdot 10^{-4}.
\]
Set $\delta =1$ and the degree of the controller $r_u$ as $3$. The data-driven controller designed by Theorem \ref{thrm:cntrldsgnCTHighorder} is
\begin{align*}
u &= x_1(  1.5x_1^2 - 0.098x_2^2  - 11.4)\\
 &\quad +x_2(0.35x_1^2  - 0.036x_2^2   - 2.0).
\end{align*}
For all $x_1,x_2\in\R$, one can bound the control input by
\[
|u| \le 2.8\sqrt 2 \|x\|^3+11\sqrt 2\|x\|.
\]

On the other hand, for all $x_1,x_2\in\R$, it holds that
\[
|R_{f_2}(x)|\le \frac{mgl}{720J}\|x\|^6 \text{ and } |R_{g_2}(x)|\le \frac{l}{6J} |x_1^3| \le  \frac{l}{6J}\|x\|^3.
\]
Then, the remainder is bounded as
\begin{align*}
|R_{f_2}(x)+R_{g_2}(x)u|&\le |R_{f_2}(x)|+|R_{g_2}(x)||u|\\
&=\! \left( \frac{mgl}{720J} + \frac{7\sqrt 2l}{15J} \right)\! \|x\|^6 \! + \frac{11\sqrt 2l}{6J} \|x\|^4\\
&\le 2.5927 ( \|x\|^6 +\|x\|^4).
\end{align*}
For the RoA estimation, we suppose that the remainder bound is over-estimated by $20\%$. More specifically, for any $x_1,x_2\in\R$, there exists a polynomial $\widetilde\rho_2(x)$ such that
\[
R_{f_2}(x)+R_{g_2}(x)u = \widetilde\rho_2(x)( \|x\|^6 +\|x\|^4)
\]
where $|\widetilde\rho_2(x)|\le 1.2\cdot2.5927=3.1112$. Let $Q_2$ be the second column of $P^{-1}$ and we obtain that $\widetilde\kappa(x)= x^{\top} Q_2( \|x\|^6 +\|x\|^4)$.

\begin{figure}
  \centering
  \includegraphics[width=0.45\textwidth]{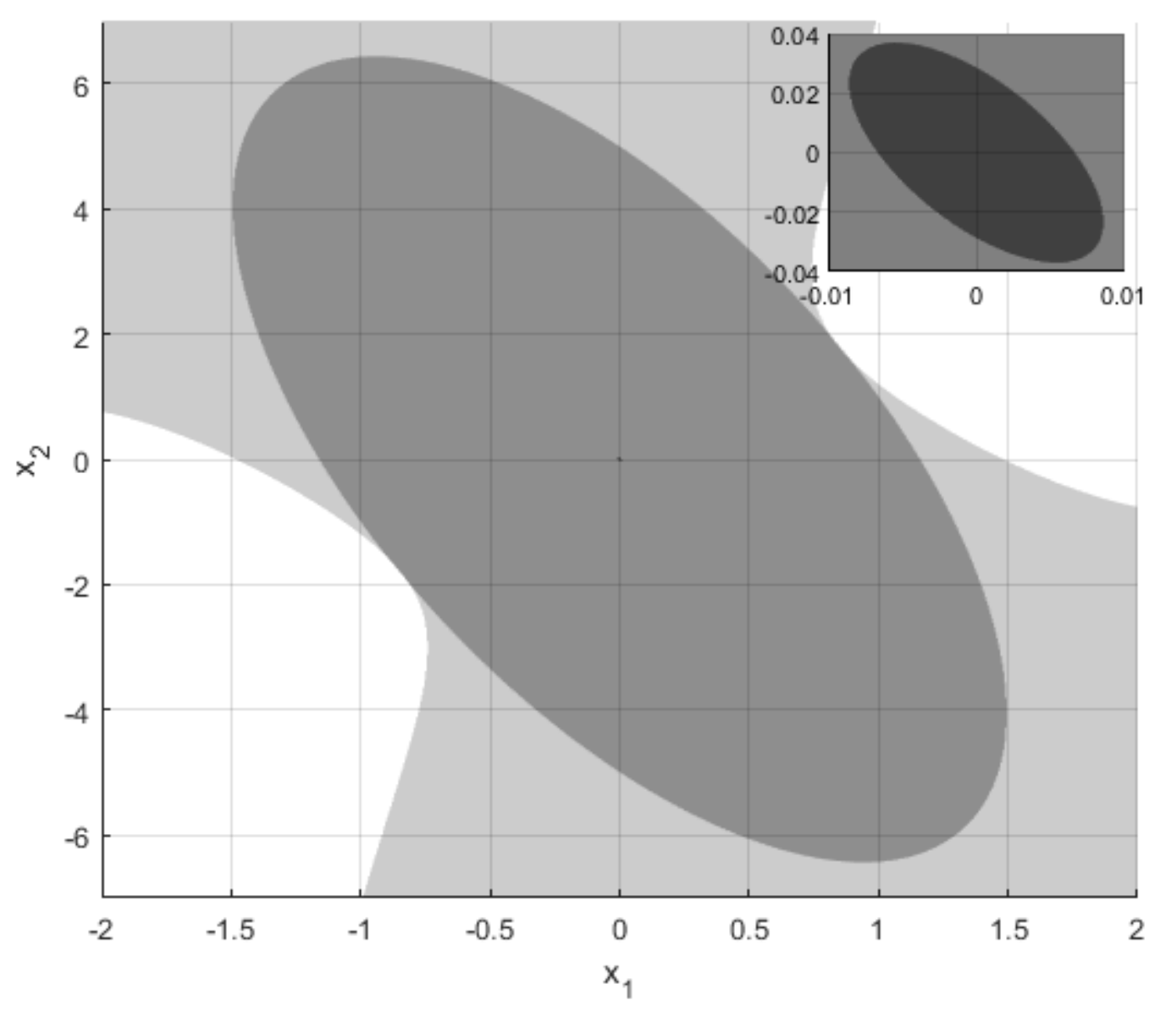}\\
  \caption{Estimations of the RoA using the Lyapunov function by the high-order approximation ($\delta=1$). The lightest grey area is the estimated RoA via the numerical method using explicit dynamics; the medium grey area is the largest sublevel set of the Lyapunov function contained in the numerically estimated RoA; the dark grey area is the estimated RoA obtained by Proposition \ref{proposition:ROAHighorder}. }\label{Fig:RoAhoAll}
\end{figure}

Using this bound and Proposition \ref{proposition:ROAHighorder}, we find an estimation of the RoA that is $\{x\in\R^2: x^{\top}P^{-1}x\le 1.58\}$ which is illustrated as the darkest area in Figure \ref{Fig:RoAhoAll}. In the same figure, the light grey area is the estimated RoA by checking point-by-point of a mesh of initial conditions using explicit dynamics, and the medium dark area is the largest sublevel set of the Lyapunov function contained in the RoA.

\begin{figure}
  \centering
  \includegraphics[width=0.45\textwidth]{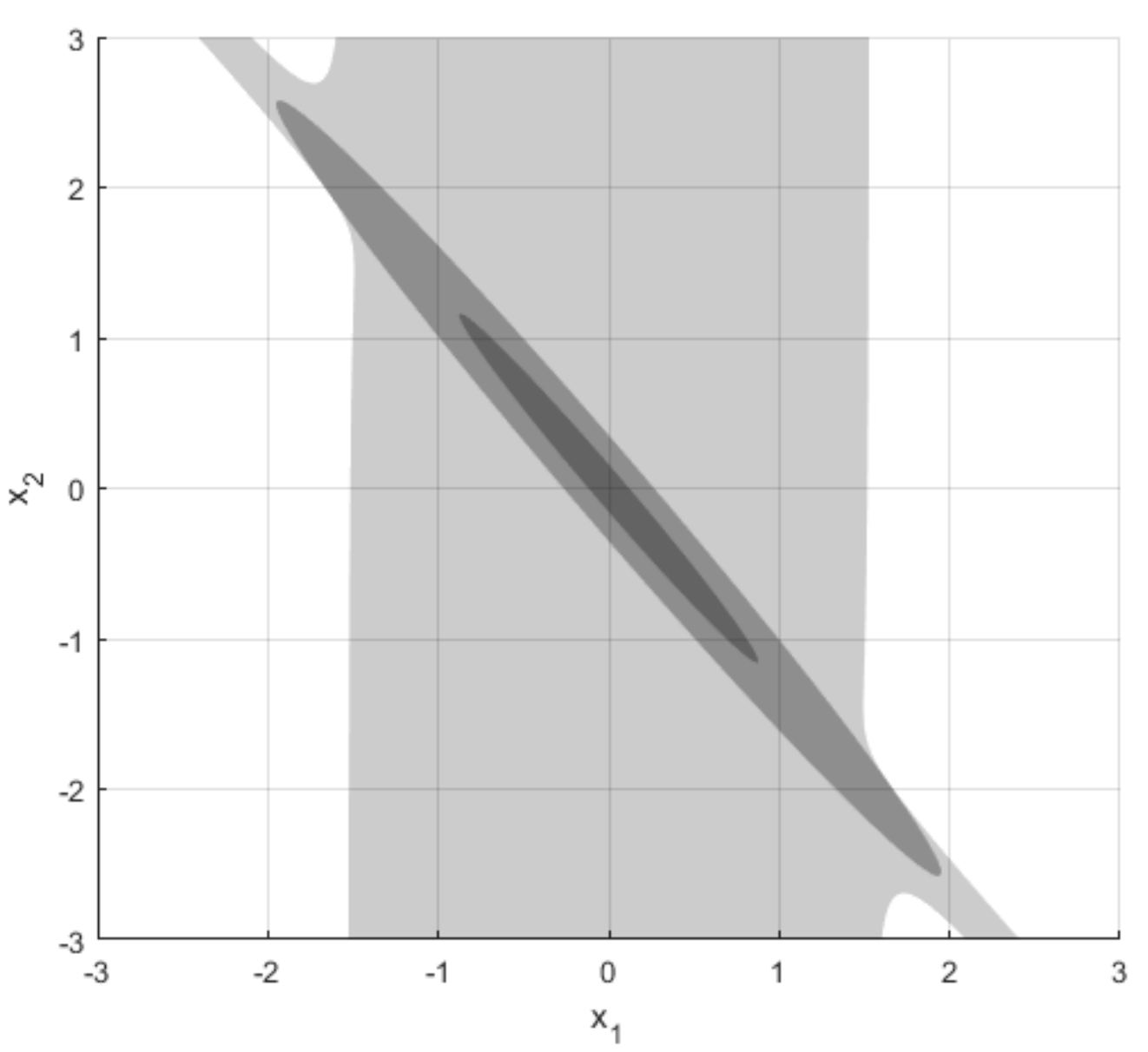}\\
  \caption{Estimations of the RoA using the Lyapunov function by the first-order approximation ($\delta=10^{-2}$). The lightest grey area is the estimated RoA using the numerical method with explicit dynamics; the medium grey area is the largest sublevel set of the Lyapunov function contained in the numerically estimated RoA; the dark grey area is the estimated RoA obtained by Proposition \ref{prpstn:ROAestimationCT}.}\label{Fig:RoA1stAll}
\end{figure}

\begin{figure}
  \centering
  \includegraphics[width=0.45\textwidth]{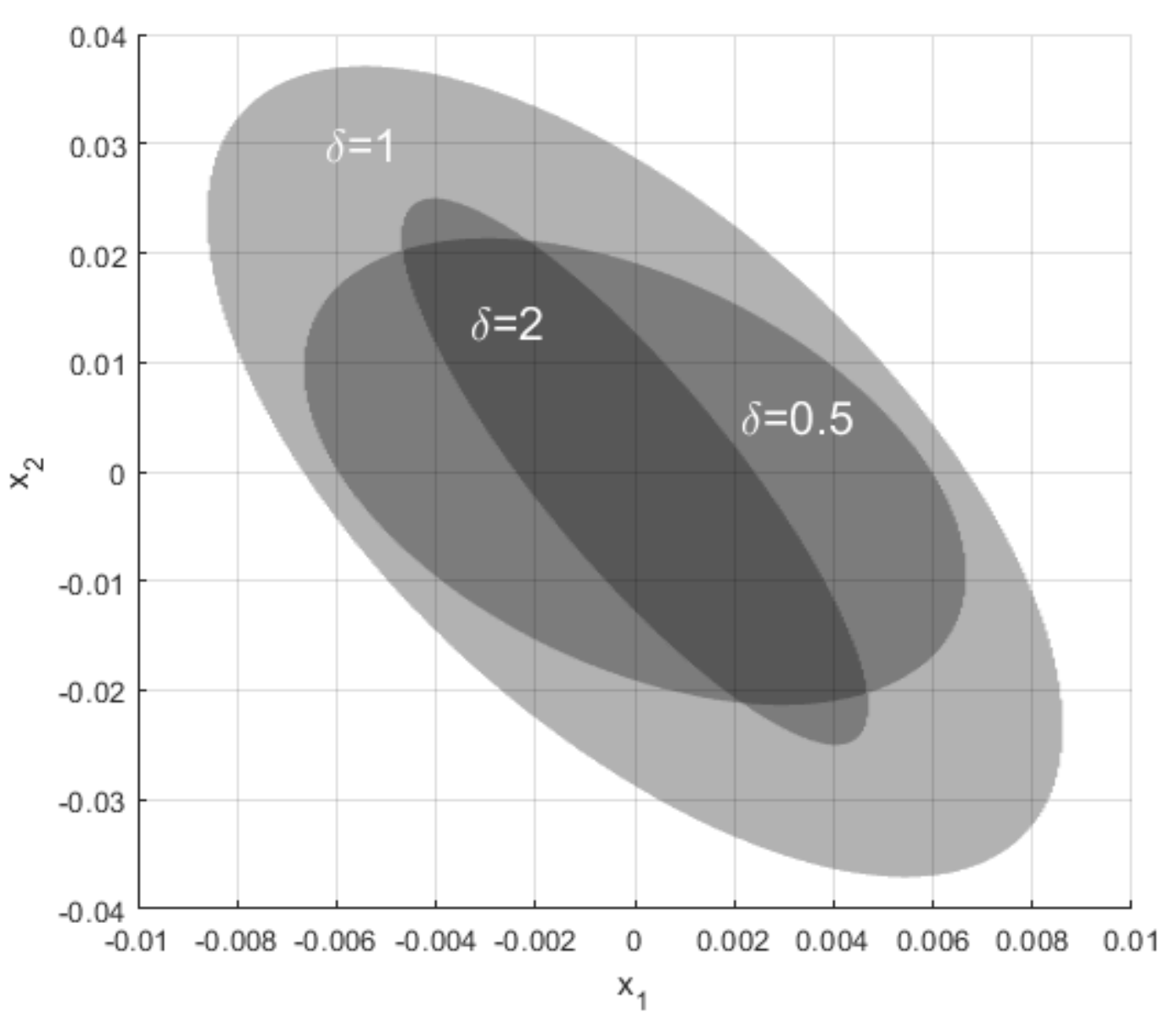}\\
  \caption{The sets $\Omega_c$ found by Proposition \ref{proposition:ROAHighorder} with different values of $\delta$.}\label{Fig:RoAHOCTdeltavary}
\end{figure}

 \begin{table*}[t]
\small
\centering
\begin{tabular}{l  c c  c}
  \toprule
   $\delta$ & $P$ & $u$ & $c^*$\\
    \toprule
  $0.5$ & $ 10^{-6}\cdot\begin{bmatrix} 0.45 & -0.63 \\-0.63 &  4.57\end{bmatrix} $ & $x_1( 0.64x_1^2 - 0.039x_2^2  - 5.15)+x_2(0.22x_1^2  - 0.027x_2^2   - 1.02)$ & $0.999$\\
  \hline
  $1$ & $ 10^{-4}\cdot \begin{bmatrix}0.47 & -1.28\\-1.28 &  8.71\end{bmatrix} $ & $x_1( 1.5x_1^2 - 0.098x_2^2  - 11.4)+x_2(0.35x_1^2  - 0.036x_2^2   - 2.0)$ & $1.58$\\
  \hline
   $ 2$  & $ 10^{-4}\cdot\begin{bmatrix}0.92 &-4.18\\ -4.18 & 25.9\end{bmatrix} $ & $x_1( 5.09x_1^2 - 0.41x_2^2  - 36.96)+x_2(0.97x_1^2  - 0.090x_2^2   - 6.33)$ &  $0.242$ \\
  \bottomrule
\end{tabular}
\caption{Simulation results on the high-order approximation with different values of $\delta$.}
 \label{table:HOvarydelta}
\end{table*}

In comparison, for the first-order approximation, the estimated RoA via the numerical method, the largest sublevel set of the Lyapunov function contained in the RoA, and the RoA estimated by Proposition \ref{prpstn:ROAestimationCT} are illustrated in Figure \ref{Fig:RoA1stAll}.

Comparing the largest RoA estimation with known dynamics and the numerical method, we see that the nonlinear controller leads to a larger RoA than the linear controller. However, if the RoA is estimated using only data, the nonlinear controller gives a much smaller area than the linear one. This could be caused by the conservativeness introduced in the process of bounding the remainder $R(x,u)$ for the high-order approximation.

For the high-order approximation, the data-driven controller design and RoA estimation are also analyzed for different values of $\delta$, and the simulation results are present in  Table \ref{table:HOvarydelta} and Figure \ref{Fig:RoAHOCTdeltavary}. In this case, we could not test values of $\delta$ in a larger range, as the SOS condition in Theorem \ref{thrm:cntrldsgnCTHighorder} is solvable by the SOSTOOLS only with $\delta$ in a small range. As shown by the simulation results, the estimations of the RoA vary with different values of $\delta$, and similar to the previous cases, how $\delta$ affects the estimation still remains as an open problem.

It should be pointed out that, besides $\delta$, there are other parameters that may affect the solving of the SOS conditions in Propositions \ref{prpstn:ROAestimationCT} and \ref{proposition:ROAHighorder}, for example, the degree of the SOS polynomials, which could also affect the RoA estimation. At the moment, it remains to be understood what is the optimal choices of the parameters that will lead to the largest estimation of the RoA. Hence, we expect that the RoA estimation for the high-order approximation case can be enlarged by carefully choosing the parameters. In this work, we aim at showing that one can use data to estimate the RoA of a nonlinear systems under a data-driven controller. Later works may look into how to enlarge the RoA estimation in the data-based setting.

\section{Conclusion}\label{section:conclusion}
For general nonlinear dynamics without explicit information on the nonlinearities, this paper proposes data-driven stabilizer designs and RoA analysis by approximating the unknown functions using Taylor's expansion. Using finite-length input-state data, linear and nonlinear stabilizers are designed for continuous-time nonlinear systems that render the known equilibrium locally asymptotically stable. Then, by estimating a bound on the Taylor remainder, data-driven conditions are given to find an invariant subset of the RoA. Simulation results on the inverted pendulum show the designed data-driven controllers and the RoA estimations for both first-order and high-order approximations. The estimation of the RoA can be conservative especially for the high-order approximation, and may be further enlarged by optimizing the choice of some parameters in the design and estimation steps. Topics such as enlarging the RoA estimation and case studies on more complicated nonlinear benchmarks are all interesting directions to be considered in future works.

\bibliographystyle{elsarticle-num}
\bibliography{referenceROA}

\appendix

%

\section{Proof of Lemma \ref{lemma:remainderbound}}\label{proofrmndrbnd}

We prove this lemma for two cases, \emph{i.e.}, $r=0$ and $r\ge 1$.

When $r=0$, Taylor's expansion of $\phi(z)$ at $z=0$ is
\begin{align*}
\phi(z)=\phi(0) + R_0(z).
\end{align*}
Then, \eqref{continuouscondition:scalar} gives $|R_0(z)|\le C\|z\|$.

For any integer $r\ge 1$, one can write the function $\phi(z)$ as
\begin{align*}
\phi(z) &= \sum_{|\alpha|\le r} \frac{\partial^{\alpha}\phi(0)}{\alpha!}z^{\alpha} +R_{r}(z) ~\text{  or}\\
\phi(z) & = \sum_{|\alpha|\le r-1} \frac{\partial^{\alpha}\phi(0)}{\alpha!}z^{\alpha} +R_{r-1}(z).
\end{align*}
As a consequence, one has that
\begin{align*}
&\quad~ R_r(z) \\
&= R_{r-1}(z) - \sum_{|\alpha|=r}\frac{\partial^{\alpha}\phi(0)}{\alpha!}z^{\alpha}\\
& = r\sum_{|\alpha|=r}\frac{z^{\alpha}}{\alpha !}\int^{1}_{0} (1-t)^{r-1} \partial^{\alpha}\phi(tz) dt- \! \sum_{|\alpha|=r}\frac{z^{\alpha}}{\alpha!} \partial^{\alpha}\phi(0)\\
& = \sum_{|\alpha|=r}\frac{z^{\alpha}}{\alpha !}~ r \int^1_0 (1-t)^{r-1}\big( \partial^{\alpha}\phi(tz) - \partial^{\alpha}\phi(0) \big) dt.
\end{align*}
By the condition \eqref{continuouscondition:scalar}, one has
\begin{align*}
|\partial^{\alpha}\phi(tz) - \partial^{\alpha}\phi(0)| \le C\|tz\|,~~t\in(0,1).
\end{align*}
Then, it holds that
\begin{align*}
\left| R_r(z) \right| &\le \sum_{|\alpha|=r}\frac{|z^{\alpha}|}{\alpha !}~ r \int^1_0 (1-t)^{r-1}Ct\|z\| dt \\
& = rC\|z\| \sum_{|\alpha|=r}\frac{|z^{\alpha}|}{\alpha !} \int^1_0 (1-t)^{r-1}tdt.
\end{align*}
Using integration by parts, one can show that
\begin{align*}
&\quad \int^1_0 (1-t)^{r-1}tdt\\
&= \frac{1}{2}\int^1_0 (1-t)^{r-1}d (t^{2})\\
&=\frac{1}{2}\left[ \left. (1-t)^{r-1}t^{2} \right|^1_0- \int^1_0 t^{2} d\big( (1-t)^{r-1}\big) \right] \\
&=\frac{r-1}{2}\int^1_0 (1-t)^{r-2}t^{2}dt\\
&~~\vdots\\
&=\frac{(r-1)!}{2 \cdots  r}\int^1_0 (1-t)^{r-r}t^r dt \\
&= \frac{(r-1)!}{(r+1)!}.
\end{align*}
By the multinomial theorem, \emph{i.e.},
\begin{align*}
\sum_{|\alpha|=r} \frac{r!}{\alpha !}z^{\alpha} = (z_1+\cdots+z_{\sigma})^{r},
\end{align*}
and the fact that $|z_1+\cdots z_{\sigma}|\le\sqrt{\sigma}\|z\|$, it holds that
\begin{align*}
\left| R_r(z) \right| &\le rC\|z\| \sum_{|\alpha|=r}\frac{|z^{\alpha}|}{\alpha !} \cdot \frac{(r-1)!}{(r+1)!}\\
&=\frac{C\|z\|}{(r+1)!}\sum_{|\alpha|=r}\frac{r!}{\alpha !}|z^{\alpha}|\\
&=\frac{C\|z\|}{(r+1)!}|z_1+\cdots z_{\sigma}|^{r}\\
&\le \frac{\sigma^{r/2} C\|z\|^{r+1}}{(r+1)!}.
\end{align*}
The proof is complete.\qed

\section{Proof of Lemma \ref{lemma:Psatz} }\label{proofPsatz}
The set inclusion condition \eqref{setcondtion:RoA} can be equivalently written as
\begin{align*}
\{ x\in\R^n: \varphi_1(x)\ge0, \varphi_2(x)\ge 0, x\ne 0\} = \emptyset.
\end{align*}
By Theorem \ref{theorem:P-satz}, we know that this is true if and only if there exist $\varphi(x)\in\mathcal S_C(\varphi_1,\varphi_2)$ and $\zeta(x)\in\mathcal S_M(x)$, such that
\begin{align}\label{equation:P-satz}
\varphi(x) + \zeta(x)^2=0.
\end{align}

Let
\[
\varphi=s_0+s_1\varphi_1+s_2\varphi_2
\]
where $s_j$, $j=0,1,2$ are SOS polynomials. By the definition of the cone $\mathcal S_C$, one has that $\varphi\in\mathcal S_C(\varphi_1,\varphi_2)$. Choosing $\zeta(x)^2=x^{\t}x$, we write the condition \eqref{equation:P-satz} as
\begin{align}\label{equation:P-satz_linear}
s_0+s_1\varphi_1+s_2\varphi_2 +x^{\t}x=0
\end{align}
As $s_0=-(s_1\varphi_1+s_2\varphi_2 +x^{\t}x)$ from \eqref{equation:P-satz_linear}, if there exist SOS polynomials $s_1$ and $s_2$ such that the SOS condition \eqref{soscondition:Psatzlemma} holds, then there exist SOS polynomials $s_j$, $j=0,1,2$ such that \eqref{equation:P-satz_linear} is true, and hence the set inclusion condition \eqref{setcondtion:RoA} holds. \qed

\section{Proof of Lemma \ref{lemma:dVCT1st}}
\label{proofdVCT1st}
For the closed-loop system with the controller $u=Kx$ designed via Theorem \ref{thrm:controldesignCT}, the derivative of the Lyapunov function $V(x)=x^{\top}P^{-1}x$ satisfies
\begin{align*}
\dot V(x) &= x^{\top}P^{-1}(A+BK)x + x^{\top}(A+BK)^{\top}P^{-1}x\\
&\quad +2x^{\top}P^{-1}R(x,Kx) \\
&\le  -wx^{\top}P^{-1}x +2x^{\top}P^{-1}R(x,Kx).
\end{align*}
Under Assumption \ref{assumption:Lipschitzr1}, for all $x\in\mathbb D$ and $i=1,\dots,n$, the bounds of the remainder can be found as
\begin{align*}
|R_i(x,Kx)| \le \frac{\sqrt{m+n} L_i}{2} \|(x,Kx)\|^{2}.
\end{align*}
%
Hence, for all $x\in\mathbb D$, there exists a continuous $\rho_i(x)$ for each $i=1,\dots,n$ such that
\begin{align*}
R_i(x,Kx) &= \rho_i(x) \|(x,Kx)\|^{2}, \\
\rho_i(x)&\in \left[ -\frac{\sqrt{m+n} L_i}{2}, \frac{\sqrt{m+n} L_i}{2}\right].
\end{align*}
Define $\rho(x) = [\rho_1(x)~\dots~\rho_n(x)]^{\top}$. By the definition of polytopes \cite[Definition 3.21]{Blanchini2008SetControl}, the vector $\rho(x)$ belongs to the polytope
\begin{align*}
\mathcal H= \{ \varrho:  -\bar h \preceq \varrho \preceq \bar h \}
\end{align*}
where
\[
\bar h=[\bar h_1~ \cdots ~\bar h_n]^{\top} = \begin{bmatrix}
 \frac{\sqrt{m+n} L_1}{2} &\cdots & \frac{\sqrt{m+n} L_n}{2}
\end{bmatrix}^{\top}.
\]
Denote $Q_i$ as the $i$th column of $P^{-1}$. It holds that
\begin{align*}
&\quad 2x^{\top}P^{-1}R(x,Kx) \\
&= 2\begin{bmatrix} x^{\top}Q_1 & \cdots & x^{\top}Q_n \end{bmatrix}
\begin{bmatrix} \rho_1(x)\|(x,Kx)\|^{2} \\ \vdots \\ \rho_n(x)\|(x,Kx)\|^{2} \end{bmatrix}\\
&= 2\sum^n_{i=1} x^{\top}Q_i\rho_i(x)\|(x,Kx)\|^{2}\\
&= 2\sum^n_{i=1} x^{\top}Q_i\|(x,Kx)\|^{2}\cdot \rho_i(x)\\
& = 2\begin{bmatrix} x^{\top}Q_1\|(x,Kx)\|^{2} & \cdots & x^{\top}Q_n\|(x,Kx)\|^{2} \end{bmatrix}\rho(x).
\end{align*}
Denote
\begin{align*}
\kappa(x)=\begin{bmatrix} x^{\top}Q_1\|(x,Kx)\|^{2} & \cdots & x^{\top}Q_n\|(x,Kx)\|^{2} \end{bmatrix}.
\end{align*}
Then, the derivative of the Lyapunov function satisfies for all $x\in\mathbb D$
\begin{align*}
\dot V(x) \le -wx^{\top}P^{-1}x + 2\kappa(x)\rho(x)
\end{align*}
where $\rho(x)\in\mathcal H$. \qed

\end{document}